\documentclass[secthm,seceqn,amsthm,ussrhead,12pt]{amsart}
\usepackage[utf8]{inputenc}
\usepackage[english]{babel}

\usepackage{amssymb,amsmath,amsthm,amsfonts,xcolor,enumerate,hyperref,comment,longtable,cleveref}
\usepackage{soul}
\usepackage{times}
\usepackage{cite}
\usepackage{pdflscape}
\usepackage{ulem}
\usepackage[mathcal]{euscript}
\usepackage{tikz}
\usepackage{hyperref}
\usepackage{cancel}
\usepackage{stmaryrd}

\usetikzlibrary{arrows}

\newtheorem{theorem}{Theorem}
\newtheorem{lemma}[theorem]{Lemma}

\newtheorem{remark}[theorem]{Remark}

\newtheorem*{theoremA}{Theorem A}
\newtheorem*{theoremB}{Theorem B}

\usepackage{stmaryrd}
\usepackage{xcolor}

\begin{document}
\noindent{\bf
One-generated  nilpotent bicommutative algebras}
\footnote{
This work was supported by  
MTM2016-79661-P;
FPU scholarship (Spain).
}

   \

   {\bf  
   Ivan   Kaygorodov$^{a}$,
   Pilar P\'{a}ez-Guill\'{a}n$^{b}$ \&
   Vasily Voronin$^{c}$}

\

{\tiny

$^{a}$ CMA, Universidade da Beira Interior, Covilhã, Portugal.

$^{b}$ University of Santiago de Compostela, Santiago de Compostela, Spain.

$^{c}$ Novosibirsk State University, Novosibirsk, Russia.

\smallskip

   E-mail addresses:

\smallskip
Ivan   Kaygorodov (kaygorodov.ivan@gmail.com)

Pilar P\'{a}ez-Guill\'{a}n (pilar.paez@usc.es)

Vasily Voronin (voronin.vasily@gmail.com)

}

\

\noindent{\bf Abstract}:
{\it We give a classification of  $5$- and $6$-dimensional complex one-generated nilpotent bicommutative algebras.}

\noindent {\bf Keywords}:
{\it bicommutative algebras, nilpotent algebras, algebraic classification, central extension.}

\noindent {\bf MSC2020}: 17A30.

\section*{Introduction}

The algebraic classification, up to isomorphism, of algebras of dimension $n$ from a certain variety
defined by a family of polynomial identities is a classical problem in the theory of non-associative algebras.
There are many results related to the algebraic classification of low-dimensional algebras in the varieties of
Jordan, Lie, Leibniz, Zinbiel and many others \cite{ack, cfk19, gkk, gkks, degr3, usefi1, degr1, degr2, ha16, hac18, kv16}.
Another interesting direction is the study of one-generated objects.
The description of one-generated, or cyclic, groups is well known: there exists a unique one-generated group of order $n$, up to isomorphism.
    In the case of algebras, similar results have been obtained in varieties such as associative \cite{karel}, non-commutative Jordan \cite{jkk19},  Leibniz and Zinbiel \cite{bakhrom}, since it was proved that there is only one $n$-dimensional one-generated nilpotent algebra from each of them.
    However, this circumstance does not hold for every variety of non-associative algebras. For example, the classifications of $4$-dimensional right commutative \cite{kkk18}, assosymmetric \cite{ikm19},  bicommutative \cite{kpv19}, commutative \cite{fkkv19} and terminal \cite{kkp19geo} nilpotent algebras show that there exist several  one-generated algebras of dimension $4$ from these varieties.
    Recently, 
    one-generated nilpotent Novikov and assosymmetric algebras in dimensions $5$ and $6,$ and one-generated nilpotent terminal algebras in dimension $5$ were classified in \cite{ckkk19, kks19, gkks}.
In the present paper, we give the algebraic classification of
$5$- and $6$-dimensional complex one-generated nilpotent bicommutative  algebras.

The variety of bicommutative algebras was introduced by Dzhumadildaev and Tulenbaev in 2003 \cite{dt03}, and it is defined by the following identities of right- and left-commutativity:
\[
\begin{array}{rclllrcl}
(xy)z &=& (xz)y, & \ &  x(yz) &=& y(xz).
\end{array} \]
One-sided commutative algebras first appeared in the paper \cite{cayley}, by Cayley, in 1857. The structure of the free bicommutative algebra of countable rank and its main numerical invariants were described by
Dzhumadildaev, Ismailov, and Tulenbaev \cite{dit11}, see also the announcement \cite{dt03}.
Bicommutative algebras were also studied in \cite{drensky1, drensky2, di}, and in \cite{DKS09, DKV10} under the name of LR-algebras.
The algebraic and geometric classifications of $2$-dimensional bicommutative algebras are given in \cite{kv16}, and that of $4$-dimensional nilpotent bicommutative algebras, in \cite{kpv19}.

Our method for classifying nilpotent bicommutative algebras is based on the calculation of central extensions of smaller nilpotent algebras from the same variety.
The algebraic study of central extensions of Lie and non-Lie algebras has been an important topic for years \cite{omirov,ha17,hac16,kkl18,ss78,krs19,zusmanovich}, especially due to their usefulness to classify nilpotent algebras from different varieties.
First, Skjelbred and Sund employed central extensions of Lie algebras to obtain a classification of nilpotent Lie algebras  \cite{ss78}.
After that, and basing on the method developed by Skjelbred and Sund,  all  non-Lie central extensions of  all $4$-dimensional Malcev algebras were described \cite{hac16},
all  anticommutative central extensions of  $3$-dimensional anticommutative algebras \cite{cfk182},
and all  central extensions of  $2$-dimensional algebras \cite{cfk18}.
Also, central extensions were used to give classifications of
  $4$-dimensional nilpotent associative algebras \cite{degr1},
  $5$-dimensional nilpotent Jordan algebras \cite{ha16},
  $5$-dimensional nilpotent restricted Lie algebras \cite{usefi1},
  $5$-dimensional nilpotent associative commutative algebras \cite{kpv19},
  $6$-dimensional nilpotent Lie algebras \cite{degr3,degr2},
  $6$-dimensional nilpotent Malcev algebras \cite{hac18},
  $6$-dimensional nilpotent anticommutative algebras \cite{kkl19},
  $8$-dimensional dual Mock Lie algebras \cite{lisa},
and some others.

The paper is organized as follows. Section~\ref{s1} is devoted to summarise the method for classifying nilpotent bicommutative algebras, inspired by that of Skjelbred and Sund~\cite{ss78}, including also some considerations about one-generated algebras and the current classification of one-generated nilpotent bicommutative algebras of dimension up to $4$. In Section~\ref{s2}, we apply this method to obtain the complete classification of $5$-dimensional one-generated nilpotent bicommutative algebras. Finally, in Section~\ref{s3}, we base on the results of the previous section to classify the one-generated nilpotent bicommutative algebras of dimension $6$.

\

\section{The algebraic classification of nilpotent bicommutative algebras}\label{s1}

\subsection{Method of classification of nilpotent algebras}
In this section, we offer an analogue of the Skjelbred-Sund method for classifying nilpotent bicommutative algebras. As other analogues of this method were carefully explained in, for instance, \cite{ hac16,cfk18}, we will limit to expose its general lines, and refer the interested reader to the previous sources. 

Let $({\bf A}, \cdot)$ be a bicommutative  algebra of dimension $n$ over  $\mathbb C$
and $\bf V$ a vector space of dimension $s$ over ${\mathbb C}$. 
Consider the bilinear maps $\theta  \colon {\bf A} \times {\bf A} \longrightarrow {\bf V}$
such that
\[ \theta(xy,z)=\theta(xz,y), \ \  \theta(x,yz)= \theta(y,xz), \]
which will be called  cocycles; they form a $\mathbb C$-linear space which we will denote by ${\rm Z}^2\left(
\bf A,\bf V \right) $.
We define also the linear subspace ${\rm B}^2\left(
{\bf A},{\bf V}\right)=\left\{\delta f\ : f\in \textup{Hom}\left( {\bf A},{\bf V}\right) \right\}$ formed by the cocycles $\delta f$ defined by  $\delta f  (x,y ) =f(xy )$. Its elements are called
 coboundaries. The  second cohomology space ${\rm H}^2\left( {\bf A},{\bf V}\right) $ is defined to be the quotient space ${\rm Z}^2
\left( {\bf A},{\bf V}\right) \big/{\rm B}^2\left( {\bf A},{\bf V}\right)$.

Let $\phi \in {\rm Aut}({\bf A}) $ be an automorphism of ${\bf A} $. Every $\theta \in
{\rm Z}^2\left( {\bf A},{\bf V}\right) $ defines another cocycle $\phi\theta$ by $\phi \theta (x,y)
=\theta \left( \phi \left( x\right) ,\phi \left( y\right) \right) $. This construction induces a right action of ${\rm Aut}({\bf A})$
 on ${\rm Z}^2\left( {\bf A},{\bf V}\right) $, which leaves
 ${\rm B}^2\left( {\bf A},{\bf V}\right) $ invariant. 
 So, it follows that ${\rm Aut}({\bf A})$ also acts on ${\rm H}^2\left( {\bf A},{\bf V}\right)$.

Given a bilinear map $\theta  \colon {\bf A} \times {\bf A} \longrightarrow {\bf V}$, we can construct the direct sum ${\bf A}_{\theta } = {\bf A}\oplus {\bf V}$ and endow it with the
bilinear product $\left[ -,-\right] _{{\bf A}_{\theta }}$ defined by $\left[ x+x^{\prime },y+y^{\prime }\right] _{{\bf A}_{\theta }}=
 xy +\theta(x,y) $ for all $x,y\in {\bf A},x^{\prime },y^{\prime }\in {\bf V}$.
It is clear that the algebra ${\bf A_{\theta}}$ is bicommutative if and only if $\theta \in {\rm Z}^2({\bf A}, {\bf V})$. Note that it is a $s$- dimensional central extension of ${\bf A}$ by ${\bf V}$; it is not difficult to prove that every central extension of $\bf A$ is of this type.

An important object for the development of the method is the so-called  annihilator of $\theta $, namely
${\rm Ann}(\theta)=\left\{ x\in {\bf A}:\theta \left( x, {\bf A} \right)+ \theta \left({\bf A} ,x\right) =0\right\} $. We also recall that the  annihilator of an  algebra ${\bf A}$ is 
the ideal ${\rm Ann}(  {\bf A} ) =\left\{ x\in {\bf A}:  x{\bf A}+ {\bf A}x =0\right\}$. It holds
 that
${\rm Ann}\left( {\bf A}_{\theta }\right) =\big({\rm Ann}(\theta) \cap{\rm Ann}({\bf A})\big)
 \oplus {\bf V}$. If ${\bf A}={\bf A}_0 \oplus I$ for any subspace $I$ of ${\rm Ann}({\bf A})$,
then $I$ is called an annihilator component of ${\bf A}$. A central extension of an algebra $\bf A$ without annihilator component is called a  non-split central extension.

The following result is fundamental for the classification method.

\begin{lemma}\label{l1}
Let ${\bf A}$ be an $n$-dimensional bicommutative algebra such that $\dim({\rm Ann}({\bf A}))=s\neq0$. Then there exists, up to isomorphism, a unique $(n-s)$-dimensional bicommutative  algebra ${\bf A}'$ and a bilinear map $\theta \in {\rm Z}^2({\bf A}, {\bf V})$ with ${\rm Ann}({\bf A'})\cap{\rm Ann}(\theta)=0$, where $\bf V$ is a vector space of dimension $s$, such that ${\bf A} \cong {{\bf A}'}_{\theta}$ and
 ${\bf A}/{\rm Ann}({\bf A})\cong {\bf A}'$.
\end{lemma}

In view of this lemma, we can solve the isomorphism problem for bicommutative algebras with non-zero annihilator just working with central extensions.

Let us fix a basis $\{e_{1},\ldots ,e_{s}\}$ of ${\bf V}$. Given a cocycle $\theta$, it can be uniquely
written as $\theta \left( x,y\right) =
\displaystyle \sum_{i=1}^{s} \theta _{i}\left( x,y\right) e_{i}$, where $\theta _{i}\in
{\rm Z}^2\left( {\bf A},\mathbb C\right) $. The main relations between $\theta$ and   $\theta_i$ are that $\theta \in
{\rm B}^2\left( {\bf A},{\bf V}\right)$ if and only if all $\theta _{i}\in {\rm B}^2\left( {\bf A},
\mathbb C\right) $, and that ${\rm Ann}(\theta)={\rm Ann}(\theta _{1})\cap\ldots \cap{\rm Ann}(\theta _{s})$. 
Furthermore, we have the following lemma.

\begin{lemma}\label{l:ind}
With the previous notations, if ${\rm Ann}(\theta)\cap {\rm Ann}\left( {\bf A}\right) =0$, then ${\bf A}_{\theta }$ has an annihilator component if and only if $\left[ \theta _{1}\right],\ldots ,\left[ \theta _{s}\right] $ are linearly dependent in ${\rm H}^2\left( {\bf A},\mathbb C\right)$.
\end{lemma}

Recall that, given a finite-dimensional vector space ${\bf V}$ over $\mathbb C$, the  Grassmannian $G_{k}\left( {\bf V}\right) $ is the set of all $k$-dimensional
linear subspaces of $ {\bf V}$. 
 Given $W=\left\langle
\left[ \theta _{1}\right],\dots,\left[ \theta _{s}
\right] \right\rangle \in G_{s}\left( {\rm H}^2\left( {\bf A},\mathbb C
\right) \right) $ and $\phi \in {\rm Aut}({\bf A})$, we define $\phi W$ as $\phi W=\left\langle \left[ \phi \theta _{1}\right],\dots,\left[ \phi \theta _{s}\right]
\right\rangle $, which again belongs to $G_{s}\left( {\rm H}^2\left( {\bf A},\mathbb C \right) \right) $. This induces an action of ${\rm Aut}({\bf A})$ on $G_{s}\left( {\rm H}^2\left( {\bf A},\mathbb C\right) \right) $; the orbit of $W\in G_{s}\left(
{\rm H}^2\left( {\bf A},\mathbb C\right) \right) $ under this action will be denoted by ${\rm Orb}(W)$. 

Consider the set 
\[
T_{s}({\bf A}) =\left\{ W=\left\langle \left[ \theta _{1}\right],\dots,\left[ \theta _{s}\right] \right\rangle \in
G_{s}\left( {\rm H}^2\left( {\bf A},\mathbb C\right) \right) : \bigcap\limits_{i=1}^{s}{\rm Ann}(\theta _{i})\cap{\rm Ann}({\bf A}) =0\right\};
\]
it is well defined, as whenever $\left\langle \left[ \theta _{1}\right],\dots,
\left[ \theta _{s}\right] \right\rangle =\left\langle \left[ \vartheta
_{1}\right],\dots,\left[ \vartheta _{s}\right]
\right\rangle \in G_{s}\left( {\rm H}^2\left( {\bf A},\mathbb C\right)
\right)$, it holds that 
\[ \bigcap\limits_{i=1}^{s}{\rm Ann}(\theta _{i})\cap {\rm Ann}\left( {\bf A}\right) = \bigcap\limits_{i=1}^{s}
{\rm Ann}(\vartheta _{i})\cap{\rm Ann}( {\bf A}).\]
In
addition, $T_{s}({\bf A})$ is stable under the action of ${\rm Aut}({\bf A})$.

Let us denote by
${\rm E}\left( {\bf A},{\bf V}\right) $ the set of all non-split $s$-dimensional central extensions of ${\bf A}$ by a $s$-dimensional linear space
${\bf V}$:
\[
{\rm E}\left( {\bf A},{\bf V}\right) =\left\{ {\bf A}_{\theta }:\theta \left( x,y\right) = \sum_{i=1}^{s}\theta _{i}\left( x,y\right) e_{i} \ \ \text{and} \ \ \left\langle \left[ \theta _{1}\right],\dots,
\left[ \theta _{s}\right] \right\rangle \in T_{s}({\bf A}) \right\} .
\]

The following lemma solves the isomorphism problem for non-split central extensions of bicommutative algebras.

\begin{lemma}\label{l2}
 Let ${\bf A}_{\theta },{\bf A}_{\vartheta }\in{\rm E}\left( {\bf A},{\bf V}\right) $. Suppose that $\theta \left( x,y\right) =  \displaystyle \sum_{i=1}^{s}
\theta _{i}\left( x,y\right) e_{i}$ and $\vartheta \left( x,y\right) =
\displaystyle \sum_{i=1}^{s} \vartheta _{i}\left( x,y\right) e_{i}$.
Then the bicommutative algebras ${\bf A}_{\theta }$ and ${\bf A}_{\vartheta } $ are isomorphic
if and only if
$${\rm Orb}\left\langle \left[ \theta _{1}\right],\dots,\left[ \theta _{s}\right] \right\rangle =
{\rm Orb}\left\langle \left[ \vartheta _{1}\right],\dots,\left[ \vartheta _{s}\right] \right\rangle .$$
\end{lemma}

In conclusion, we can construct all the non-split central extensions of a bicommutative algebra ${\bf A'}$ of
dimension $n-s$ by following this procedure:

\begin{enumerate}
\item Determine ${\rm H}^2( {\bf A}',\mathbb {C}) $, ${\rm Ann}({\bf A}')$ and ${\rm Aut}({\bf A}')$.

\item Determine the set of ${\rm Aut}({\bf A}')$-orbits on $T_{s}({\bf A}') $.

\item For each orbit, construct the bicommutative algebra associated with a
representative of it.
\end{enumerate}

Taking into account Lemmas~\ref{l1} and~\ref{l2}, it is clear that, provided that the classifications of nilpotent bicommutative algebras of dimension up to $n-1$ are known, we can also classify all the nilpotent bicommutative algebras of dimension $n$.

Note also that if we want to stick to the one-generated case, it suffices to consider the non-split central extensions of one-generated nilpotent bicommutative algebras of lower dimension. Indeed, the central extensions of an algebra which is not one-generated cannot be one-generated; on the other hand, considering the definition of ${\rm B}^2\left({\bf A},{\bf V}\right)$ and Lemma~\ref{l:ind}, it is not difficult to see that the non-split extensions of a one-generated algebra are again one-generated.

\subsection{Notations}
Let ${\bf A}$ be a bicommutative algebra and let us fix 
a basis $\{e_{1},\dots,e_{n}\}$. We will denote by $\Delta _{ij}$ the
bicommutative bilinear form
$\Delta _{ij} \colon {\bf A} \times {\bf A}\longrightarrow \mathbb C$
with $\Delta _{ij}\left( e_{l},e_{m}\right) = \delta_{il}\delta_{jm}$.
Then, the set $\left\{ \Delta_{ij}:1\leq i, j\leq n\right\} $ is a basis for the linear space of
the bilinear forms on ${\bf A}$, and every $\theta \in
{\rm Z}^2\left( {\bf A},\bf V\right)$ can be uniquely written as $
\theta = \displaystyle \sum_{1\leq i,j\leq n} c_{ij}\Delta _{{i}{j}}$, where $
c_{ij}\in \mathbb C$.

Henceforth, ${\mathcal B}^i_j$ will denote the $j^{\textup{th}}$ $i$-dimensional one-generated bicommutative algebra.

\subsection{Low dimensional one-generated nilpotent bicommutative algebras}
We can extract from \cite{kpv19} the list of the one-generated nilpotent bicommutative algebras of dimension up to $4$. Note that, for dimensions $3$ and $4$, it is not necessary to check the lists of ``trivial'' algebras (i.e. those in which the identity $(xy)z=x(yz)=0$ holds), since they cannot be one-generated.

\[\begin{array}{ll llllllllllll}
{\mathcal B}^2_{01} &:& e_1 e_1 = e_2 \\
{\mathcal B}^3_{01} &:& e_1 e_1 = e_2 & e_2 e_1 = e_3 \\
{\mathcal B}^3_{02}(\lambda) &:& e_1 e_1 = e_2 & e_1 e_2 = e_3 & e_2 e_1 = \lambda e_3 \\
{\mathcal B}^4_{01} &:& e_1 e_1 = e_2 & e_1 e_2 = e_4 & e_2 e_1 = e_3 \\
{\mathcal B}^4_{02} &:& e_1 e_1 = e_2 & e_1 e_2 = e_4 & e_2 e_1 = e_3 & e_3 e_1 = e_4  \\
{\mathcal B}^4_{03} &:& e_1 e_1 = e_2 & e_2 e_1 = e_3 & e_3 e_1 = e_4  \\
{\mathcal B}^4_{04} &:& e_1 e_1 = e_2 & e_1 e_2 = e_3 & e_1 e_3 = e_4 & e_2 e_1 = e_4  \\
{\mathcal B}^4_{05} &:& e_1 e_1 = e_2 & e_1 e_2 = e_3 & e_1 e_3 = e_4 & e_2 e_1 = e_3 + e_4 & e_2 e_2 = e_4 & e_3 e_1 = e_4  \\
{\mathcal B}^4_{06}(\lambda) &:& e_1 e_1 = e_2 & e_1 e_2=e_3 & e_1 e_3=e_4 & e_2 e_1 = \lambda e_3 & e_2 e_2=\lambda e_4 & e_3 e_1=\lambda e_4. \\
\end{array} \]

\

\section{Classification of $5$-dimensional one-generated nilpotent bicommutative algebras}\label{s2}

There are not $3$-dimensional central extensions of the $2$-dimensional one-generated bicommutative algebra ${\mathcal B}^2_{01}$.

\subsection{$2$-dimensional central extensions of $3$-dimensional one-generated algebras}
Considering 2-dimensional central extensions of 3-dimensional one-generated nilpotent bicommutative algebras ${\mathcal B}^3_{01}$ and ${\mathcal B}^3_{02}(\lambda)$ we get the algebras ${\mathcal B}^5_{01}$ and ${\mathcal B}^5_{02}(\lambda)$:
\begin{longtable}{ll lllllllll} ${\mathcal B}^5_{01}$ &$:$ & $e_1 e_1 = e_2$  & $e_1e_2=e_4$ & $e_2e_1=e_3$ & $e_3e_1=e_5$ & \\
 ${\mathcal B}^5_{02}(\lambda)$ &$:$& $e_1 e_1 = e_2$ & $e_1 e_2=e_3$ & $e_1e_3=e_5$ & $e_2 e_1=\lambda e_3+e_4$ &  $e_2e_2=\lambda e_5$ & $e_3e_1=\lambda e_5$. \\
\end{longtable}

\subsection{Cohomology spaces  of $4$-dimensional one-generated bicommutative algebras}
In the following table we give the description of the second cohomology space of  $4$-dimensional one-generated nilpotent bicommutative algebras.

\

\

\begin{longtable}{|lll|}

\hline

${\rm Z^2}({\mathcal B}^4_{01})$ &=& $\Big\langle
\begin{array}{l}\Delta_{11},\Delta_{12}, \Delta_{13}+\Delta_{22}+\Delta_{41}, \Delta_{14}, \Delta_{21},\Delta_{31}  \end{array}\Big\rangle$ \\

${\rm B^2}({\mathcal B}^4_{01})$ &=& $\Big\langle \begin{array}{l}\Delta_{11}, \Delta_{12}, \Delta_{21} \end{array}\Big\rangle$ \\

${\rm H^2}({\mathcal B}^4_{01})$ &=& $\Big\langle\begin{array}{l}[\Delta_{14}], [\Delta_{13}] + [\Delta_{22}]+ [\Delta_{41}], [\Delta_{31}] \end{array}\Big\rangle$   \\

\hline

${\rm Z^2}({\mathcal B}^4_{02})$ &=& $\Big\langle  \begin{array}{l}\Delta_{11}, \Delta_{12}, \Delta_{13}+\Delta_{22}+\Delta_{41}, \Delta_{21}, \Delta_{31}  \end{array}\Big\rangle$ \\

${\rm B^2}({\mathcal B}^4_{02})$ &=& $\Big\langle \begin{array}{l}\Delta_{11}, \Delta_{12}+\Delta_{31}, \Delta_{21} \end{array}\Big\rangle$ \\

${\rm H^2}({\mathcal B}^4_{02})$ &=& $\Big\langle \begin{array}{l} [\Delta_{13}] + [\Delta_{22}] + [\Delta_{41}], [\Delta_{31}] \end{array}\Big\rangle$   \\

\hline

${\rm Z^2}({\mathcal B}^4_{03})$ &=& $\Big\langle \begin{array}{l} \Delta_{11},\Delta_{12}, \Delta_{21}, \Delta_{31}, \Delta_{41} \end{array}\Big\rangle$ \\

${\rm B^2}({\mathcal B}^4_{03})$ &=& $\Big\langle \begin{array}{l}\Delta_{11}, \Delta_{21}, \Delta_{31} \end{array}\Big\rangle$ \\

${\rm H^2}({\mathcal B}^4_{03})$ &=& $\Big\langle \begin{array}{l}[\Delta_{12}], [\Delta_{41}] \end{array}\Big\rangle$   \\

\hline

${\rm Z^2}({\mathcal B}^4_{04})$ &=& $\Big\langle  \begin{array}{l}\Delta_{11}, \Delta_{12}, \Delta_{13}, \Delta_{14}+\Delta_{22}+\Delta_{31}, \Delta_{21} \end{array}\Big\rangle$ \\

${\rm B^2}({\mathcal B}^4_{04})$ &=& $\Big\langle \begin{array}{l}\Delta_{11}, \Delta_{12}, \Delta_{13}+\Delta_{21} \end{array}\Big\rangle$ \\

${\rm H^2}({\mathcal B}^4_{04})$ &=& $\Big\langle \begin{array}{l} [\Delta_{14}]+[\Delta_{22}]+[\Delta_{31}], [\Delta_{21}] \end{array}\Big\rangle$   \\

\hline

${\rm Z^2}({\mathcal B}^4_{05})$ &=& $\Big\langle \begin{array}{l} \Delta_{11}, \Delta_{12}, \Delta_{13}+\Delta_{22}+\Delta_{31},  \Delta_{14}+\Delta_{22}+\Delta_{23}+\Delta_{31}+\Delta_{32}+\Delta_{41}, \Delta_{21} \end{array}\Big\rangle$ \\

${\rm B^2}({\mathcal B}^4_{05})$ &=& $\Big\langle \begin{array}{l} \Delta_{11}, \Delta_{12}+\Delta_{21}, \Delta_{13}+\Delta_{21}+\Delta_{22}+\Delta_{31} \end{array}\Big\rangle$ \\

${\rm H^2}({\mathcal B}^4_{05})$ &=& $\Big\langle \begin{array}{l} [\Delta_{14}]+[\Delta_{22}]+[\Delta_{23}]+[\Delta_{31}]+[\Delta_{32}]+[\Delta_{41}], [\Delta_{21}] \end{array}\Big\rangle$   \\

\hline

${\rm Z^2}({\mathcal B}^4_{06}(\lambda))$ &=& $\Big\langle \begin{array}{l} \Delta_{11},\Delta_{12}, \Delta_{13}+\lambda\Delta_{22}+\lambda\Delta_{31}, \Delta_{14}+\lambda\Delta_{23}+\lambda\Delta_{32}+\lambda\Delta_{41}, \Delta_{21} \end{array}\Big\rangle$ \\

${\rm B^2}({\mathcal B}^4_{06}(\lambda))$ &=& $\Big\langle \begin{array}{l} \Delta_{11}, \Delta_{12}+\lambda\Delta_{21}, \Delta_{13}+\lambda\Delta_{22}+\lambda\Delta_{31} \end{array}\Big\rangle$ \\

${\rm H^2}({\mathcal B}^4_{06}(\lambda))$ &=& $\Big\langle \begin{array}{l} [\Delta_{14}]+\lambda[\Delta_{23}]+\lambda[\Delta_{32}]+\lambda[\Delta_{41}], [\Delta_{21}] \end{array}\Big\rangle$   \\

\hline

\end{longtable}

\subsection{Central extensions of ${\mathcal B}^4_{01}$}

Let us use the following notations:
\[\nabla_1=[\Delta_{14}], \nabla_2=[\Delta_{31}], \nabla_3=[\Delta_{13}]+[\Delta_{41}]+[\Delta_{22}].\]

The automorphism group of ${\mathcal B}^4_{01}$ consists of invertible matrices of the form

\[\phi=\left(
                             \begin{array}{cccc}
                               x & 0  & 0 & 0 \\
                               y & x^2  & 0 & 0 \\
                               z & xy & x^3 & 0 \\
                               t & xy & 0 & x^3                           \end{array}\right)
                               .\]
                               
Since
\[
\phi^T
                           \left(\begin{array}{cccc}
                                0 & 0 & \alpha_3 & \alpha_1  \\
                                0 & \alpha_3 & 0 & 0  \\
                                \alpha_2 & 0 & 0 & 0 \\
                                \alpha_3 & 0 & 0 & 0
                             \end{array}
                           \right)\phi
                           =\left(\begin{array}{cccc}
                                 \alpha^* & \alpha^{**} & \alpha_3^* & \alpha_1^*   \\
                                 \alpha^{***} & \alpha_3^* & 0 & 0  \\
                                 \alpha_2^* & 0 & 0 & 0 \\
                                 \alpha_3^* & 0 & 0 & 0
                             \end{array}\right),\]
we have that the action of ${\rm Aut} ({\mathcal B}^4_{01})$ on the subspace
$\langle  \sum\limits_{i=1}^3\alpha_i \nabla_i \rangle$
is given by
$\langle  \sum\limits_{i=1}^3\alpha^*_i \nabla_i \rangle,$ where

\[
\alpha^*_1=x^4\alpha_1,\quad \alpha^*_2=x^4\alpha_2,\quad \alpha^*_3=x^4\alpha_3.
\]

\subsubsection{$1$-dimensional central extensions}
We have the following cases:
\begin{enumerate}
\item if $\alpha_3\neq 0$, by choosing $x=\frac{1}{\sqrt[4]{\alpha_3}}$ we have the representative $\langle \lambda \nabla_1 +\mu \nabla_2 +  \nabla_3 \rangle;$
\item if $\alpha_3=0$, then it must hold that $\alpha_1\neq 0$ (otherwise, the new algebras would have $2$-dimensional annihilator and could be constructed as $2$-dimensional central extensions of $3$-dimensional bicommutative algebras), and by choosing $x=\frac{1}{\sqrt[4]{\alpha_1}}$, we have the representatives 
$\langle \nabla_1+\lambda\nabla_2\rangle$.
\end{enumerate}

Hence, we obtain  the following new algebras: 

\begin{longtable}{lllllll lllll} ${\mathcal B}^5_{03}(\lambda, \mu)$ &:&
$e_1 e_1 = e_2$ & $e_1 e_2=e_4$ & $e_1 e_3= e_5$ & $e_1 e_4= \lambda e_5$ \\
&&  $e_2e_1=e_3$ & $e_2 e_2= e_5$ & $e_3 e_1=\mu e_5$  & $e_4e_1= e_5$&& \\
 ${\mathcal B}^5_{04}(\lambda)$ &:&
$e_1 e_1 = e_2$ & $e_1 e_2=e_4$ & $e_1 e_4= e_5$ &  $e_2e_1= e_3$ & $e_3 e_1=\lambda e_5$. \\ 
\end{longtable}

\subsubsection{$2$-dimensional central extensions}
Consider the vector space generated by the following two cocycles:
\[\begin{array}{rcl}
\theta_1 &=& \alpha_1 \nabla_1+\alpha_2\nabla_2+\alpha_3\nabla_3  \\
\theta_2 &=& \beta_1 \nabla_1+\beta_2\nabla_2.
\end{array}\]

If $\alpha_3=0$, we get the representative $\langle \nabla_1, \nabla_2 \rangle.$
If $\alpha_3\neq 0$, we distinguish the following cases:

\begin{enumerate}
    \item if $\beta_2=0$,  we have the family of representatives 
    $\langle \nabla_1, \lambda \nabla_2+\nabla_3 \rangle;$
    
    \item if $\beta_2\neq 0$,  we have the family of representatives 
    $\langle \lambda \nabla_1+\nabla_2, \mu \nabla_1 + \nabla_3 \rangle.$
\end{enumerate}

Hence, we have the following new algebras:

\begin{longtable}{ll lllllll} ${\mathcal B}^6_{01}$ &:&
$e_1 e_1 = e_2$ & $e_1e_2=e_4$ & $e_1 e_4=e_5$ & $e_2 e_1=e_3$ & $e_3 e_1=e_6$ \\
 ${\mathcal B}^6_{02}(\lambda)$ &:&
$e_1 e_1 = e_2$ &   $e_1 e_2=e_4$ & $e_1e_3=e_6$ &  $e_1 e_4=e_5$ \\
&&  $e_2 e_1= e_3$ & $e_2e_2= e_6$  & $e_3e_1=\lambda e_6$ & $e_4e_1= e_6$  \\
 ${\mathcal B}^6_{03}(\lambda, \mu)$ &:& 
$e_1 e_1 = e_2$ &   $e_1 e_2=e_4$ &  $e_1e_3=  e_6$ &   \multicolumn{2}{l}{$e_1 e_4=\lambda e_5 + \mu e_6$} \\
&& $e_2 e_1= e_3$  & $e_2e_2= e_6$  & $e_3e_1= e_5$ & $e_4e_1= e_6$.  &&
\end{longtable}

\subsection{Central extensions of ${\mathcal B}^4_{02}$}

Let us use the following notations:
\[\nabla_1=[\Delta_{31}], \nabla_2=[\Delta_{13}]+[\Delta_{22}]+[\Delta_{41}].\]

The automorphism group of ${\mathcal B}^4_{02}$ consists of invertible matrices of the form

\[\phi=\left(
                             \begin{array}{cccc}
                               1 & 0  & 0 & 0 \\
                               x & 1  & 0 & 0 \\
                               y & x & 1 & 0 \\
                               z & x+y & x & 1                           \end{array}\right)
                               .\]
                               
Since
\[
\phi^T
                           \left(\begin{array}{cccc}
                                0 & 0 & \alpha_2 & 0  \\
                                0 & \alpha_2 & 0 & 0  \\
                                \alpha_1 & 0 & 0 & 0 \\
                                \alpha_2 & 0 & 0 & 0
                             \end{array}
                           \right)\phi
                           =\left(\begin{array}{cccc}
                                 \alpha^* & \alpha^{**} & \alpha_2^* & 0 \\
                                 \alpha^{***} & \alpha_2^* & 0 & 0  \\
                                 \alpha_1^*+\alpha^{**} & 0 & 0 & 0 \\
                                 \alpha_2^* & 0 & 0 & 0
                             \end{array}\right),\]
we have that the action of ${\rm Aut} ({\mathcal B}^4_{02})$ on the subspace
$\langle  \sum\limits_{i=1}^2\alpha_i \nabla_i \rangle$
is given by
$\langle  \sum\limits_{i=1}^2\alpha^*_i \nabla_i \rangle,$ where

\[
\alpha^*_1=\alpha_1-x\alpha_2,\quad \alpha^*_2=\alpha_2.
\]

\subsubsection{$1$-dimensional central extensions}
We can suppose that $ \alpha_2\neq0$, as otherwise we would obtain algebras with $2$-dimensional annihilator which would be $2$-dimensional central extensions of $3$-dimensional bicommutative algebras. Then, by choosing $x=\frac{\alpha_1}{\alpha_2}$, we have the representative $\langle \nabla_2 \rangle$, whose associated algebra is

\begin{longtable}{lllll lllll}$ {\mathcal B}^5_{05}$ &:&
$e_1 e_1 = e_2$  & $e_1e_2=e_4$ & $e_1 e_3=e_5$ & $e_2 e_1=e_3$ & $e_2e_2=e_5$ & $e_3 e_1=e_4$ &  $e_4 e_1=e_5$. \\
\end{longtable}

\subsubsection{$2$-dimensional central extensions}
Since the cohomology space has dimension 2, we obtain the following new algebra: 

\begin{longtable}{ll llllllll} ${\mathcal B}^6_{04}$ &:&
$e_1 e_1 = e_2$ &    $e_1 e_2=e_4$ & $e_1 e_3=e_6$ & $e_2e_1=e_3$ & $e_2e_2=e_6$ & $e_3e_1=e_4+e_5$ & $e_4 e_1=e_6$. \\
\end{longtable}

\subsection{Central extensions of ${\mathcal B}^4_{03}$}

Let us use the following notations:
\[\nabla_1=[\Delta_{12}], \nabla_2=[\Delta_{41}].\]

The automorphism group of ${\mathcal B}^4_{03}$ consists of invertible matrices of the form

\[\phi=\left(
                             \begin{array}{cccc}
                               x & 0  & 0 & 0 \\
                               y & x^2  & 0 & 0 \\
                               z & xy & x^3 & 0 \\
                               t & xz & x^2y & x^4                           \end{array}\right)
                               .\]
                               
Since
\[
\phi^T
                           \left(\begin{array}{cccc}
                                0 & \alpha_1 & 0 & 0  \\
                                0 & 0 & 0 & 0  \\
                                0 & 0 & 0 & 0 \\
                                \alpha_2 & 0 & 0 & 0
                             \end{array}
                           \right)\phi
                           =\left(\begin{array}{cccc}
                                 \alpha^* & \alpha_1^* & 0 & 0 \\
                                 \alpha^{**} & 0 & 0 & 0  \\
                                 \alpha^{***} & 0 & 0 & 0 \\
                                 \alpha_2^* & 0 & 0 & 0
                             \end{array}\right),\]
we have that the action of ${\rm Aut} ({\mathcal B}^4_{03})$ on the subspace
$\langle  \sum\limits_{i=1}^2\alpha_i \nabla_i \rangle$
is given by
$\langle  \sum\limits_{i=1}^2\alpha^*_i \nabla_i \rangle,$ where

\[
\alpha^*_1=x^3\alpha_1 \quad \alpha^*_2=x^5\alpha_2.
\]

\subsubsection{$1$-dimensional central extensions}
We can suppose again that $ \alpha_2\neq0$. So, we have the following cases:
\begin{enumerate}
\item if $\alpha_1\neq 0$, by choosing $x=\sqrt{\frac{\alpha_1}{\alpha_2}}$ we have the representative $\langle \nabla_1+\nabla_2 \rangle$;
\item if $\alpha_1=0$,  by choosing $x=\frac{1}{\sqrt[3]{\alpha_2}}$ we have the representative $\langle \nabla_2 \rangle$.
\end{enumerate}

Hence, we get the following new algebras: 

\begin{longtable}{ll lllll} ${\mathcal B}^5_{06}$ &:&
$e_1 e_1 = e_2$ & $e_1e_2=e_5$ & $e_2 e_1=e_3$ & $e_3 e_1=e_4$ & $e_4 e_1= e_5$\\ 
 ${\mathcal B}^5_{07}$ &:&
$e_1 e_1 = e_2$ &   $e_2 e_1=e_3$ &    $e_3 e_1=e_4$ & $e_4 e_1=e_5$. \\
\end{longtable}

\subsubsection{$2$-dimensional central extensions}
Since the cohomology space has dimension 2, we have the following new algebra: 

\begin{longtable}{ll lllll} ${\mathcal B}^6_{05}$ &:&
$e_1 e_1 = e_2$ & $e_1e_2=e_5$ &   $e_2 e_1=e_3$ &    $e_3 e_1=e_4$  & $e_4 e_1=e_6$. \\ 
\end{longtable}

\subsection{Central extensions of ${\mathcal B}^4_{04}$}

Let us use the following notations:
\[\nabla_1=[\Delta_{21}], \nabla_2=[\Delta_{14}]+[\Delta_{22}]+[\Delta_{31}].\]

The automorphism group of ${\mathcal B}^4_{04}$ consists of invertible matrices of the form

\[\phi=\left(
                             \begin{array}{cccc}
                               1 & 0  & 0 & 0 \\
                               x & 1  & 0 & 0 \\
                               y & x & 1 & 0 \\
                               z & x+y & x & 1                           \end{array}\right)
                               .\]
                               
Since
\[
\phi^T
                           \left(\begin{array}{cccc}
                                0 & 0 & 0 & \alpha_2  \\
                                \alpha_1 & \alpha_2 & 0 & 0  \\
                                \alpha_2 & 0 & 0 & 0 \\
                                0 & 0 & 0 & 0
                             \end{array}
                           \right)\phi
                           =\left(\begin{array}{cccc}
                                 \alpha^* & \alpha^{**} & \alpha^{***} & \alpha_2^* \\
                                 \alpha_1^*+\alpha^{***} & \alpha_2^* & 0 & 0  \\
                                 \alpha_2^* & 0 & 0 & 0 \\
                                 0 & 0 & 0 & 0
                             \end{array}\right),\]
we have that the action of ${\rm Aut} ({\mathcal B}^4_{04})$ on the subspace
$\langle  \sum\limits_{i=1}^2\alpha_i \nabla_i \rangle$
is given by
$\langle  \sum\limits_{i=1}^2\alpha^*_i \nabla_i \rangle,$ where
\[
\alpha^*_1=\alpha_1+x\alpha_2,\quad \alpha^*_2=\alpha_2.
\]

\subsubsection{$1$-dimensional central extensions}
We can suppose that $ \alpha_2\neq0$. Then, by choosing $x=-{\frac{\alpha_1}{\alpha_2}}$, we have the representative $\langle \nabla_2 \rangle$, with the following associated algebra:

\begin{longtable}{ll lllllll} ${\mathcal B}^5_{08}$ &:& 
$e_1 e_1 = e_2$ &   $e_1 e_2=e_3$ & $e_1 e_3=e_4$ & $e_1e_4=e_5$ & $e_2e_1=e_4$ & $e_2 e_2=e_5$ & $e_3e_1=e_5$. \\
\end{longtable}

\subsubsection{$2$-dimensional central extensions}
Since the cohomology space has dimension 2, we get  the following new algebra: 

\begin{longtable}{llll lllll} ${\mathcal B}^6_{06}$ &:& 
$e_1 e_1 = e_2$ &  $e_1 e_2=e_3$ & $e_1 e_3=e_4$ & $e_1e_4=e_6$ & $e_2e_1=e_4+e_5$ & $e_2e_2=e_6$ & $e_3e_1=e_6$. \\
\end{longtable}

\subsection{Central extensions of ${\mathcal B}^4_{05}$}

Let us use the following notations:
\[\nabla_1=[\Delta_{21}], \nabla_2=[\Delta_{14}]+[\Delta_{22}]+[\Delta_{23}]+[\Delta_{31}]+[\Delta_{32}]+[\Delta_{41}].\]

The automorphism group of ${\mathcal B}^4_{05}$ consists of invertible matrices of the form

\[\phi=\left(
                             \begin{array}{cccc}
                               1 & 0  & 0 & 0 \\
                               x & 1  & 0 & 0 \\
                               y & 2x & 1 & 0 \\
                               z & x^2+x+2y & 3x & 1                           \end{array}\right)
                               .\]
                               
Since
\[
\phi^T
                           \left(\begin{array}{cccc}
                                0 & 0 & 0 & \alpha_2  \\
                                \alpha_1 & \alpha_2 & \alpha_2 & 0  \\
                                \alpha_2 & \alpha_2 & 0 & 0 \\
                                \alpha_2 & 0 & 0 & 0
                             \end{array}
                           \right)\phi
                           =\left(\begin{array}{cccc}
                                 \alpha^* & \alpha^{**} & \alpha^{***} & \alpha_2^* \\
                                 \alpha_1^*+\alpha^{**}+\alpha^{***} & \alpha_2^*+\alpha^{***} & 0 & 0  \\
                                 \alpha_2^*+\alpha^{***} & 0 & 0 & 0 \\
                                 \alpha_2^* & 0 & 0 & 0
                             \end{array}\right),\]
we have that the action of ${\rm Aut} ({\mathcal B}^4_{05})$ on the subspace
$\langle  \sum\limits_{i=1}^2\alpha_i \nabla_i \rangle$
is given by
$\langle  \sum\limits_{i=1}^2\alpha^*_i \nabla_i \rangle,$ where

\[
\alpha^*_1=\alpha_1-2x\alpha_2,\quad \alpha^*_2=\alpha_2.
\]

\subsubsection{$1$-dimensional central extensions}
We can suppose that $ \alpha_2\neq0$. Then, by choosing $x={\frac{\alpha_1}{2\alpha_2}}$, we have the representative $\langle \nabla_2 \rangle$. Hence, we obtain  the following new algebra: 

\begin{longtable}{ll lllllll} ${\mathcal B}^5_{09}$ &:& 
$e_1 e_1 = e_2$&   $e_1 e_2=e_3$ & $e_1 e_3=e_4$ & $e_1e_4=e_5$ & $e_2e_1=e_3+e_4$  \\
&& $e_2e_2=e_4+e_5$   & $e_2e_3=e_5$ & $e_3e_1=e_4+e_5$ & $e_3e_2=e_5$& $e_4e_1=e_5$.\\
\end{longtable}

\subsubsection{$2$-dimensional central extensions}
Since the cohomology space has dimension 2, we have  the following new algebra: 

\begin{longtable}{ll lllllll} ${\mathcal B}^6_{07}$ &:&
$e_1 e_1 = e_2$ &   $e_1 e_2=e_3$ & $e_1 e_3=e_4$  & $e_1e_4=e_6$ & $e_2e_1= e_3+e_4+e_5$ \\
&& $e_2e_2=e_4+e_6$   & $e_2e_3=e_6$ & $e_3e_1=e_4+e_6$  & $e_3e_2=e_6$ & $e_4e_1=e_6$.
\end{longtable}

\subsection{Central extensions of ${\mathcal B}^4_{06}(\lambda)$}

Let us use the following notations:
\[\nabla_1=[\Delta_{21}], \nabla_2=[\Delta_{14}] + \lambda[\Delta_{23}] + \lambda[\Delta_{32}] + \lambda[\Delta_{41}].\]

The automorphism group of ${\mathcal B}^4_{06}(\lambda)$ consists of invertible matrices of the form

\[\phi=\left(
                             \begin{array}{cccc}
                               x & 0  & 0 & 0 \\
                               0 & x^2  & 0 & 0 \\
                               y & 0 & x^3 & 0 \\
                               z & (1+\lambda)xy & 0 & x^4                           \end{array}\right)
                               .\]
                               
Since
\[
\phi^T
                           \left(\begin{array}{cccc}
                                0 & 0 & 0 & \alpha_2  \\
                                \alpha_1 & 0 & \lambda\alpha_2 & 0  \\
                                0 & \lambda\alpha_2 & 0 & 0 \\
                                \lambda\alpha_2 & 0 & 0 & 0
                             \end{array}
                           \right)\phi
                           =\left(\begin{array}{cccc}
                                 \alpha^* & \alpha^{**} & 0 & \alpha_2^* \\
                                 \alpha_1^*+\lambda\alpha^{**} & 0 & \lambda\alpha_2^* & 0  \\
                                 0 & \lambda\alpha_2^* & 0 & 0 \\
                                 \lambda\alpha_2^* & 0 & 0 & 0
                             \end{array}\right),\]
we have that the action of ${\rm Aut} ({\mathcal B}^4_{06}(\lambda))$ on the subspace
$\langle  \sum\limits_{i=1}^2\alpha_i \nabla_i \rangle$
is given by
$\langle  \sum\limits_{i=1}^2\alpha^*_i \nabla_i \rangle,$ where

\[
\alpha^*_1=x^3\alpha_1+(1-\lambda)\lambda x^2y\alpha_2,\quad \alpha^*_2=x^5\alpha_2.
\]

\subsubsection{$1$-dimensional central extensions}
We can suppose that $ \alpha_2\neq0$. So, we have the following cases:
\begin{enumerate}
\item if $\lambda \neq 0, 1$, by choosing $y=-\frac{x\alpha_1}{(1-\lambda)\lambda\alpha_2}$ and $x=\frac{1}{\sqrt[5]{\alpha_2}}$ we have the representative $\langle \nabla_2 \rangle$;
\item if $\lambda =0$ or $\lambda=1$, and $\alpha_1=0$, by choosing $x=\frac{1}{\sqrt[5]{\alpha_2}}$ we get again the representative $\langle \nabla_2 \rangle$;
\item if $\lambda =0$ or $\lambda=1$, and $\alpha_1\neq0$, by choosing $x=\sqrt{\frac{\alpha_1}{\alpha_2}}$ we have the representative $\langle \nabla_1+\nabla_2 \rangle$.
\end{enumerate}

Hence, we obtain the family of algebras ${\mathcal B}^5_{10}(\lambda)$, associated with the representative $\langle \nabla_2 \rangle$, and the algebras ${\mathcal B}^5_{11}$ and ${\mathcal B}^5_{12}$ associated with $\langle \nabla_1 + \nabla_2 \rangle$ for the values $\lambda=0$ and $\lambda=1$, respectively.

\begin{longtable}{ll llllllll} ${\mathcal B}^5_{10}(\lambda)$ &:&
$e_1 e_1 = e_2$ & $e_1 e_2=e_3$ & $e_1 e_3=e_4$ & $e_1e_4=e_5$ & $e_2e_1= \lambda e_3$ \\
&& $e_2 e_2= \lambda e_4$   & $e_2e_3=\lambda e_5$ & $e_3e_1= \lambda e_4$ & $e_3e_2=\lambda e_5$ & $e_4e_1=\lambda e_5$\\
 ${\mathcal B}^5_{11}$ &:&
$e_1 e_1 = e_2$ &   $e_1 e_2=e_3$ &    $e_1 e_3=e_4$  & $e_1e_4=e_5$ & $e_2e_1=e_5$ \\ 
 ${\mathcal B}^5_{12}$ &:&
$e_1 e_1 = e_2$ &   $e_1 e_2=e_3$ & $e_1 e_3=e_4$ & $e_1e_4=e_5$ & $e_2e_1=e_3+e_5$ \\
&& $e_2e_2=e_4$  & $e_2e_3=e_5$ & $e_3 e_1= e_4$ & $e_3e_2=e_5$ & $e_4e_1=e_5$.\\
\end{longtable}

\subsubsection{$2$-dimensional central extensions}
Since the cohomology space has dimension 2, we get  the following new algebra: 

\begin{longtable}{ll lllllll} ${\mathcal B}^6_{08}(\lambda)$ &:& 
$e_1 e_1 = e_2$ &   $e_1 e_2=e_3$ & $e_1 e_3=e_4$ & $e_1e_4=e_6$ & $e_2e_1=\lambda e_3+e_5$ \\
&& $e_2e_2=\lambda e_4$  & $e_2e_3=\lambda e_6$ &  $e_3e_1=\lambda e_4$  & $e_3e_2=\lambda e_6$ & $e_4e_1=\lambda e_6$.  \\
\end{longtable}

\ 

\subsection{Classification theorem}
The results of the previous subsections yield the following theorem.

\begin{theoremA}
Let $\mathcal B$ be a $5$-dimensional complex one-generated nilpotent bicommutative algebra.
Then, $\mathcal B$ is isomorphic to an algebra  from the following list:

\begin{longtable}{ l|lllll }

\hline
${\mathcal B}^5_{01}$ & $e_1 e_1 = e_2$  & $e_1e_2=e_4$ & $e_2e_1=e_3$ & $e_3e_1=e_5$ & \\

\hline
${\mathcal B}^5_{02}(\lambda)$ & $e_1 e_1 = e_2$ & $e_1 e_2=e_3$ & $e_1e_3=e_5$ \\& $e_2 e_1=\lambda e_3+e_4$ &  $e_2e_2=\lambda e_5$ 
& $e_3e_1=\lambda e_5$ \\

\hline
${\mathcal B}^5_{03}(\lambda, \mu)$ &
$e_1 e_1 = e_2$ & $e_1 e_2=e_4$ & $e_1 e_3= e_5$ & $e_1 e_4= \lambda e_5$ \\
&  $e_2e_1=e_3$ & $e_2 e_2= e_5$ & $e_3 e_1=\mu e_5$ & $e_4e_1= e_5$  \\

\hline
${\mathcal B}^5_{04}(\lambda)$ &
$e_1 e_1 = e_2$ & $e_1 e_2=e_4$ & $e_1 e_4= e_5$ &  $e_2e_1= e_3$ & $e_3 e_1=\lambda e_5$ \\ 

\hline
${\mathcal B}^5_{05}$ &
$e_1 e_1 = e_2$  & $e_1e_2=e_4$ & $e_1 e_3=e_5$ & $e_2 e_1=e_3$ \\
& $e_2e_2=e_5$ 
& $e_3 e_1=e_4$ &  $e_4 e_1=e_5$  \\

\hline
${\mathcal B}^5_{06}$ &
$e_1 e_1 = e_2$ & $e_1e_2=e_5$ & $e_2 e_1=e_3$ & $e_3 e_1=e_4$ & $e_4 e_1= e_5$\\ 

\hline
${\mathcal B}^5_{07}$ &
$e_1 e_1 = e_2$ &   $e_2 e_1=e_3$ &    $e_3 e_1=e_4$ & $e_4 e_1=e_5$ & \\

\hline
${\mathcal B}^5_{08}$ & 
$e_1 e_1 = e_2$ &   $e_1 e_2=e_3$ & $e_1 e_3=e_4$ & $e_1e_4=e_5$ \\
& $e_2e_1=e_4$ 
& $e_2 e_2=e_5$ & $e_3e_1=e_5$  \\

\hline
${\mathcal B}^5_{09}$ & 
$e_1 e_1 = e_2$&   $e_1 e_2=e_3$ & $e_1 e_3=e_4$ & $e_1e_4=e_5$ & $e_2e_1=e_3+e_4$ \\ 
 & $e_2e_2=e_4+e_5$ & $e_2e_3=e_5$ & $e_3e_1=e_4+e_5$ & $e_3e_2=e_5$& $e_4e_1=e_5$\\

\hline
${\mathcal B}^5_{10}(\lambda)$ &
$e_1 e_1 = e_2$ & $e_1 e_2=e_3$ & $e_1 e_3=e_4$ & $e_1e_4=e_5$ & $e_2e_1= \lambda e_3$ \\ 
& $e_2 e_2= \lambda e_4$ & $e_2e_3=\lambda e_5$ & $e_3e_1= \lambda e_4$ & $e_3e_2=\lambda e_5$ & $e_4e_1=\lambda e_5$\\

\hline
${\mathcal B}^5_{11}$ &
$e_1 e_1 = e_2$ &   $e_1 e_2=e_3$ &    $e_1 e_3=e_4$  & $e_1e_4=e_5$ & $e_2e_1=e_5$ \\ 

\hline
${\mathcal B}^5_{12}$ &
$e_1 e_1 = e_2$ &   $e_1 e_2=e_3$ & $e_1 e_3=e_4$ & $e_1e_4=e_5$ & $e_2e_1=e_3+e_5$ \\ 
& $e_2e_2=e_4$ & $e_2e_3=e_5$ & $e_3 e_1= e_4$ & $e_3e_2=e_5$ & $e_4e_1=e_5$\\

\hline

\end{longtable}
\end{theoremA}

\

\section{Classification of $6$-dimensional one-generated nilpotent bicommutative algebras}\label{s3}

\subsection{Cohomology spaces  of $5$-dimensional one-generated bicommutative algebras}
The multiplication tables of all the $5$-dimensional one-generated nilpotent bicommutative algebras are given in Theorem A (see previous section).
The necessary information about 
coboundaries, cocycles and  second cohomology spaces of these algebras is displayed in the following table.

{\tiny
\begin{longtable}{|lll|}

\hline

${\rm Z^2}({\mathcal B}^5_{01})$ &=& $\Big\langle \begin{array}{l}
\Delta_{11},\Delta_{12},\Delta_{13}+\Delta_{22}+\Delta_{41}, \Delta_{14},
\Delta_{21}, \Delta_{31}, \Delta_{51} \end{array} \Big\rangle$ \\

${\rm B^2}({\mathcal B}^5_{01})$ &=& $\Big\langle \begin{array}{l} \Delta_{11}, \Delta_{12}, \Delta_{21}, \Delta_{31} \end{array} \Big\rangle$ \\

${\rm H^2}({\mathcal B}^5_{01})$ &=& $\Big \langle \begin{array}{l} [\Delta_{13}]+[\Delta_{22}]+[\Delta_{41}], [\Delta_{14}], [\Delta_{51}]\relax \end{array} \Big\rangle$   \\

\hline

${\rm Z^2}({\mathcal B}^5_{02}(\lambda))$ &=& $\Big\langle \begin{array}{l}
\Delta_{11},\Delta_{12},\Delta_{13}-\lambda\Delta_{14}, \Delta_{14}+\Delta_{22}+\Delta_{31}, \Delta_{15}+\lambda\Delta_{23}+\lambda\Delta_{32}+\lambda\Delta_{51},
\Delta_{21}, \Delta_{41} \end{array} \Big\rangle$ \\

${\rm B^2}({\mathcal B}^5_{02}(\lambda))$ &=& $\Big\langle \begin{array}{l}\Delta_{11}, \Delta_{12}, \Delta_{13}+\lambda\Delta_{22}+\lambda\Delta_{31}, \Delta_{21} \end{array} \Big\rangle$ \\

${\rm H^2}({\mathcal B}^5_{02}(\lambda))$ &=& $ \Big \langle \begin{array}{l} [\Delta_{14}]+[\Delta_{22}]+[\Delta_{31}], [\Delta_{15}]+\lambda[\Delta_{23}]+\lambda[\Delta_{32}]+\lambda[\Delta_{51}], [\Delta_{41}] \relax \end{array} \Big \rangle$   \\

\hline

${\rm Z^2}({\mathcal B}^5_{03}(\lambda, \mu)_{\lambda=0\text{ or }\mu\neq 1/\lambda})$ &=& $\Big\langle \begin{array}{l}
\Delta_{11},\Delta_{12}, \Delta_{13}+\Delta_{22}+\Delta_{41}, \Delta_{14}, \Delta_{21}, \Delta_{31} \end{array} \Big\rangle$ \\

${\rm B^2}({\mathcal B}^5_{03}(\lambda, \mu)_{\lambda=0\text{ or }\mu\neq 1/\lambda})$ &=& $\Big \langle \begin{array}{l}\Delta_{11}, \Delta_{12},  \Delta_{13}+\lambda\Delta_{14}+\Delta_{22}+\mu\Delta_{31}+\Delta_{41}, \Delta_{21} \end{array} \Big \rangle$ \\

${\rm H^2}({\mathcal B}^5_{03}(\lambda, \mu)_{\lambda=0\text{ or }\mu\neq 1/\lambda})$ &=& $\Big \langle \begin{array}{l} [\Delta_{14}], [\Delta_{31}] \end{array} \Big \rangle$   \\

\hline

${\rm Z^2}({\mathcal B}^5_{03}(\lambda, 1/\lambda)_{\lambda\neq 0})$ &=& $\Big\langle \begin{array}{l}
\Delta_{11},\Delta_{12}, \Delta_{13}+\Delta_{22}+\Delta_{41},\Delta_{14}, \lambda \Delta_{15}+\Delta_{23}+\lambda\Delta_{24}+\Delta_{32}+\lambda\Delta_{42}+\Delta_{51}, \Delta_{21}, \Delta_{31} \end{array} \Big\rangle$ \\

${\rm B^2}({\mathcal B}^5_{03}(\lambda, 1/\lambda)_{\lambda\neq 0})$ &=& $\Big \langle \begin{array}{l}\Delta_{11}, \Delta_{12},  \Delta_{13}+\lambda\Delta_{14}+\Delta_{22}+(1/\lambda)\Delta_{31}+\Delta_{41}, \Delta_{21} \end{array} \Big \rangle$ \\

${\rm H^2}({\mathcal B}^5_{03}(\lambda, 1/\lambda)_{\lambda\neq 0})$ &=& $\Big \langle \begin{array}{l} [\Delta_{14}], \lambda [\Delta_{15}]+[\Delta_{23}]+\lambda[\Delta_{24}]+[\Delta_{32}]+\lambda[\Delta_{42}]+[\Delta_{51}], [\Delta_{31}] \end{array} \Big \rangle$   \\

\hline

${\rm Z^2}({\mathcal B}^5_{04}(\lambda))$ &=& $\Big\langle \begin{array}{l}
\Delta_{11},\Delta_{12}, \Delta_{13}+\Delta_{22}+\Delta_{41}, \Delta_{14}, \Delta_{21}, \Delta_{31} \end{array} \Big\rangle$ \\

${\rm B^2}({\mathcal B}^5_{04}(\lambda))$ &=& $\Big\langle \begin{array}{l}\Delta_{11}, \Delta_{12}, \Delta_{14}+\lambda\Delta_{31}, \Delta_{21} \end{array} \Big\rangle$ \\

${\rm H^2}({\mathcal B}^5_{04}(\lambda))$ &=& $\Big \langle \begin{array}{l} [\Delta_{13}]+[\Delta_{22}]+[\Delta_{41}],  [\Delta_{31}] \end{array} \Big\rangle$   \\

\hline

${\rm Z^2}({\mathcal B}^5_{05})$ &=& $\Big\langle \begin{array}{l}
\Delta_{11},\Delta_{12},\Delta_{13}+\Delta_{22}+\Delta_{41}, \Delta_{14}+\Delta_{23}+\Delta_{32}+\Delta_{51}, \Delta_{21}, \Delta_{31} \end{array} \Big\rangle$ \\

${\rm B^2}({\mathcal B}^5_{05})$ &=& $\Big\langle \begin{array}{l}\Delta_{11}, \Delta_{12}+\Delta_{31}, \Delta_{13}+\Delta_{22}+\Delta_{41}, \Delta_{21} \end{array} \Big\rangle$ \\

${\rm H^2}({\mathcal B}^5_{05})$ &=& $\Big \langle \begin{array}{l}   [\Delta_{14}]+[\Delta_{23}]+[\Delta_{32}]+[\Delta_{51}], [\Delta_{31}] \end{array} \Big\rangle$   \\

\hline

${\rm Z^2}({\mathcal B}^5_{06})$ &=& $\Big\langle \begin{array}{l}
\Delta_{11},\Delta_{12},\Delta_{13}+\Delta_{22}+\Delta_{51}, \Delta_{21}, \Delta_{31}, \Delta_{41} \end{array} \Big\rangle$ \\

${\rm B^2}({\mathcal B}^5_{06})$ &=& $\Big\langle \begin{array}{l} \Delta_{11}, \Delta_{12}+\Delta_{41}, \Delta_{21}, \Delta_{31} \end{array} \Big\rangle$ \\

${\rm H^2}({\mathcal B}^5_{06})$ &=& $\Big \langle \begin{array}{l}  [\Delta_{13}]+[\Delta_{22}]+[\Delta_{51}], [\Delta_{41}] \end{array} \Big\rangle$   \\

\hline

${\rm Z^2}({\mathcal B}^5_{07})$ &=& $\Big\langle \begin{array}{l}
\Delta_{11}, \Delta_{12}, \Delta_{21}, \Delta_{31}, \Delta_{41}, \Delta_{51} \end{array} \Big\rangle$ \\

${\rm B^2}({\mathcal B}^5_{07})$ &=& $\Big\langle \begin{array}{l} \Delta_{11}, \Delta_{21}, \Delta_{31}, \Delta_{41} \end{array} \Big\rangle$ \\

${\rm H^2}({\mathcal B}^5_{07})$ &=& $\Big \langle \begin{array}{l} [\Delta_{12}], [\Delta_{51}] \end{array} \Big\rangle$   \\

\hline

${\rm Z^2}({\mathcal B}^5_{08})$ &=& $\Big\langle \begin{array}{l}
\Delta_{11}, \Delta_{12}, \Delta_{13}, \Delta_{14}+\Delta_{22}+\Delta_{31}, \Delta_{15}+\Delta_{23}+\Delta_{32}+\Delta_{41}, \Delta_{21} \end{array} \Big\rangle$ \\

${\rm B^2}({\mathcal B}^5_{08})$ &=& $\Big\langle \begin{array}{l} \Delta_{11}, \Delta_{12}, \Delta_{13}+\Delta_{21}, \Delta_{14}+\Delta_{22}+\Delta_{31} \end{array} \Big\rangle$ \\

${\rm H^2}({\mathcal B}^5_{08})$ &=& $\Big \langle \begin{array}{l}  [\Delta_{15}]+[\Delta_{23}]+[\Delta_{32}]+[\Delta_{41}], [\Delta_{21}] \end{array} \Big\rangle$   \\

\hline

${\rm Z^2}({\mathcal B}^5_{09})$ &=& $\Big\langle \begin{array}{l}
\Delta_{11}, \Delta_{12}, \Delta_{13}+\Delta_{22}+\Delta_{31}, \Delta_{14}+\Delta_{22}+\Delta_{23}+\Delta_{31}+\Delta_{32}+\Delta_{41}, \\ \Delta_{15}+\Delta_{23}+\Delta_{24}+\Delta_{32}+\Delta_{33}+\Delta_{41}+\Delta_{42}+\Delta_{51}, \Delta_{21} \end{array} \Big\rangle$ \\

${\rm B^2}({\mathcal B}^5_{09})$ &=& $\Big\langle \begin{array}{l} \Delta_{11}, \Delta_{12}+\Delta_{21}, \Delta_{13}+\Delta_{21}+\Delta_{22}+\Delta_{31}, \Delta_{14}+\Delta_{22}+\Delta_{23}+\Delta_{31}+\Delta_{32}+\Delta_{41} \end{array} \Big\rangle$ \\

${\rm H^2}({\mathcal B}^5_{09})$ &=& $\Big \langle \begin{array}{l}   [\Delta_{15}]+[\Delta_{23}]+[\Delta_{24}]+[\Delta_{32}]+[\Delta_{33}]+[\Delta_{41}]+ [\Delta_{42}]+[\Delta_{51}], [\Delta_{21}] \end{array} \Big\rangle$   \\

\hline

${\rm Z^2}({\mathcal B}^5_{10}(\lambda))$ &=& $\Big\langle \begin{array}{l}
\Delta_{11}, \Delta_{12}, \Delta_{13}+\lambda\Delta_{22}+\lambda\Delta_{31}, \Delta_{14}+\lambda\Delta_{23}+\lambda\Delta_{32}+\lambda\Delta_{41}, \\ \Delta_{15}+\lambda\Delta_{24}+\lambda\Delta_{33}+\lambda\Delta_{42}+\lambda\Delta_{51}, \Delta_{21} \end{array} \Big\rangle$ \\

${\rm B^2}({\mathcal B}^5_{10}(\lambda))$ &=& $\Big\langle \begin{array}{l} \Delta_{11}, \Delta_{12}+\lambda\Delta_{21}, \Delta_{13}+\lambda\Delta_{22}+\lambda\Delta_{31}, \Delta_{14}+\lambda\Delta_{23}+\lambda\Delta_{32}+\lambda\Delta_{41} \end{array} \Big\rangle$ \\

${\rm H^2}({\mathcal B}^5_{10}(\lambda))$ &=& $\Big \langle \begin{array}{l}   [\Delta_{15}]+\lambda[\Delta_{24}]+\lambda[\Delta_{33}]+\lambda[\Delta_{42}]+\lambda[\Delta_{51}], [\Delta_{21}] \end{array} \Big\rangle$   \\

\hline

${\rm Z^2}({\mathcal B}^5_{11})$ &=& $\Big\langle \begin{array}{l}
\Delta_{11}, \Delta_{12}, \Delta_{13}, \Delta_{14}, \Delta_{15}+\Delta_{22}+\Delta_{31}, \Delta_{21} \end{array} \Big\rangle$ \\

${\rm B^2}({\mathcal B}^5_{11})$ &=& $\Big \langle \begin{array}{l} \Delta_{11}, \Delta_{12}, \Delta_{13}, \Delta_{14}+\Delta_{21} \end{array} \Big \rangle$ \\

${\rm H^2}({\mathcal B}^5_{11})$ &=& $\Big \langle \begin{array}{l}  [\Delta_{15}]+[\Delta_{22}]+[\Delta_{31}], [\Delta_{21}] \end{array} \Big \rangle$   \\

\hline

${\rm Z^2}({\mathcal B}^5_{12})$ &=& $\Big\langle \begin{array}{l}
\Delta_{11}, \Delta_{12}, \Delta_{13}+\Delta_{22}+\Delta_{31}, \Delta_{14}+\Delta_{23}+\Delta_{32}+\Delta_{41}, \\ \Delta_{15}+\Delta_{22}+\Delta_{24}+\Delta_{31}+\Delta_{33}+\Delta_{42}+\Delta_{51}, \Delta_{21} \end{array} \Big\rangle$ \\

${\rm B^2}({\mathcal B}^5_{12})$ &=& $\Big\langle \begin{array}{l} \Delta_{11}, \Delta_{12}+\Delta_{21}, \Delta_{13}+\Delta_{22}+\Delta_{31}, \Delta_{14}+\Delta_{21}+\Delta_{23}+\Delta_{32}+\Delta_{41} \end{array} \Big\rangle$ \\

${\rm H^2}({\mathcal B}^5_{12})$ &=& $\Big \langle \begin{array}{l}   [\Delta_{15}]+[\Delta_{22}]+[\Delta_{24}]+[\Delta_{31}]+[\Delta_{33}]+[\Delta_{42}]+[\Delta_{51}], [\Delta_{21}] \end{array} \Big\rangle$   \\

\hline

\end{longtable}
}

\begin{remark}
The  extensions of the algebras ${\mathcal B}^5_{03}(\lambda, \mu)_{\lambda=0\text{ or }\mu\neq 1/\lambda}$ and  ${\mathcal B}^5_{04}(\lambda)$ have $2$-dimensional annihilator, and have already been constructed as $2$-dimensional central extensions of $4$-dimensional bicommutative algebras. Then, in the following subsections we will study only the central extensions of the other algebras.
\end{remark}

\subsection{Central extensions of ${\mathcal B}^5_{01}$}

Let us use the following notations:
\[\nabla_1=[\Delta_{14}], \quad \nabla_2=[\Delta_{51}], \quad  \nabla_3=[\Delta_{13}]+[\Delta_{22}]+[\Delta_{41}].\]

The automorphism group of ${\mathcal B}^5_{01}$ consists of invertible matrices of the form
\[\phi=\begin{pmatrix}
x & 0 & 0  & 0& 0\\
y & x^2 & 0 & 0& 0\\
z & xy & x^3 & 0& 0 \\
t & xy & 0 &  x^3& 0\\
s & xz & x^2y & 0 & x^4\\
\end{pmatrix}. \]

Since
\[ \phi^T\begin{pmatrix}
0 &0 &  \alpha_3 & \alpha_1 & 0\\
0 & \alpha_3 & 0 & 0 & 0 \\
0 & 0 & 0 & 0 & 0\\
\alpha_3 & 0 & 0 & 0 & 0\\
\alpha_2 & 0 & 0 & 0 & 0
\end{pmatrix}\phi =
\begin{pmatrix}
\alpha^{*} & \alpha^{**} & \alpha_3^* &  \alpha_1^* & 0 \\
\alpha^{***} & \alpha_3^* &  0 & 0 & 0 \\
\alpha^{****} & 0 & 0 & 0 &0\\
\alpha_3^* & 0 & 0 & 0 &0\\
\alpha_2^* & 0 & 0 & 0 &0\\
\end{pmatrix},
\]
 we have that the action of ${\rm Aut} ({\mathcal B}^5_{01})$ on the subspace
$\langle \sum\limits_{i=1}^3 \alpha_i\nabla_i  \rangle$
is given by
$\langle \sum\limits_{i=1}^3 \alpha_i^* \nabla_i \rangle,$
where
\[\alpha_1^*=x^4\alpha_1, \quad \alpha_2^* =x^5\alpha_2, \quad \alpha_3^* =x^4\alpha_3. \]

We can suppose that $ \alpha_2\neq0$; otherwise, we would obtain algebras with $2$-dimensional annihilator which have already been constructed as $2$-dimensional central extensions of $4$-dimensional bicommutative algebras. We have the following cases:
\begin{enumerate}
\item if $\alpha_3 \neq 0$,  by choosing $x=\frac{\alpha_3}{\alpha_2}$ we have the family of representatives $\langle \lambda\nabla_1 + \nabla_2 + \nabla_3 \rangle$;

\item if $\alpha_3=0, \alpha_1\neq0$,  by choosing $x= \frac{\alpha_1} {\alpha_2}$ we have the representative $\langle \nabla_1 + \nabla_2 \rangle$;

\item if $\alpha_3 = 0, \alpha_1 = 0$, by choosing $x=\frac{1}{\sqrt[4]{\alpha_2}}$ we have the representative $\langle \nabla_2 \rangle$.

\end{enumerate}

Hence, we obtain the following new algebras:

\begin{longtable}{ll lllllllllll} ${\mathcal B}^6_{09}(\lambda)$ &:& 
$e_1 e_1 = e_2$  & $e_1 e_2=e_4$ & $e_1e_3=e_6$ & $e_1e_4=\lambda e_6$ & $e_2 e_1=e_3$ \\
&& $e_2e_2=e_6$ & $e_3e_1=e_5$  & $e_4e_1=e_6$ & $e_5e_1=e_6$ \\
 ${\mathcal B}^6_{10}$ &:& 
$e_1 e_1 = e_2$ & $e_1 e_2=e_4$ & $e_1e_4=e_6$ & $e_2 e_1=e_3$ & $e_3e_1=e_5$ & $e_5e_1=e_6$ \\ 
 ${\mathcal B}^6_{11}$ &:& 
$e_1 e_1 = e_2$  & $e_1 e_2=e_4$ & $e_2 e_1=e_3$ & $e_3e_1=e_5$ & $e_5e_1=e_6$. \\ 
\end{longtable}

\subsection{Central extensions of ${\mathcal B}^5_{02}(\lambda)$}

Let us use the following notations:
\[\nabla_1 = [\Delta_{41}], \quad \nabla_2 = [\Delta_{14}] + [\Delta_{22}] + [\Delta_{31}], \quad \nabla_3 = [\Delta_{15}] + \lambda[\Delta_{23}] + \lambda[\Delta_{32}] + \lambda[\Delta_{51}].\]

The automorphism group of ${\mathcal B}^5_{01}$ consists of invertible matrices of the form
\[\phi=\begin{pmatrix}
x & 0 & 0  & 0 & 0\\
y & x^2 & 0 & 0 & 0\\
z & (1+\lambda)xy & x^3 & 0 & 0 \\
t & xy & 0 & x^3 & 0\\
s & \lambda y^2+(1+\lambda)xz  & (1+2\lambda)x^2y & \lambda(1-\lambda)x^2y & x^4\\
\end{pmatrix}. \]

Since
\[ \phi^T\begin{pmatrix}
0 & 0 &  0 & \alpha_2 & \alpha_3\\
0 & \alpha_2 & \lambda\alpha_3 & 0 & 0 \\
\alpha_2 & \lambda\alpha_3 & 0 & 0 & 0\\
\alpha_1 & 0 & 0 & 0 & 0\\
\lambda\alpha_3 & 0 & 0 & 0 & 0
\end{pmatrix}\phi =
\begin{pmatrix}
\alpha^{*} & \alpha^{**} & \alpha^{***} & \alpha_2^* & \alpha_3^* \\
\alpha^{****} & \alpha_2^*+\lambda \alpha^{***} & \lambda\alpha_3^* & 0 & 0 \\
\alpha_2^*+\lambda \alpha^{***} & \lambda\alpha_3^* & 0 & 0 &0\\
\alpha_1^{*} & 0 & 0 & 0 &0\\
\lambda\alpha_3^* & 0 & 0 & 0 &0\\
\end{pmatrix},
\]
 we have that the action of ${\rm Aut} ({\mathcal B}^5_{02}(\lambda))$ on the subspace
$\langle \sum\limits_{i=1}^3 \alpha_i\nabla_i  \rangle$
is given by
$\langle \sum\limits_{i=1}^3 \alpha_i^* \nabla_i \rangle,$
where
\[\alpha_1^*=x^4\alpha_1 + (1-\lambda)\lambda^2x^3y\alpha_3, \quad \alpha_2^* =x^4\alpha_2 + (1-\lambda)\lambda x^3y\alpha_3, \quad \alpha_3^* = x^5\alpha_3. \]

We can suppose that $ \alpha_3\neq0$. We have the following cases:
\begin{enumerate}
\item if $\lambda = 0$ or $\lambda = 1$, and  $\alpha_1 =\alpha_2= 0$, then we have the representative $\langle \nabla_3 \rangle;$

\item if $\lambda = 0$ or $\lambda = 1$,  $\alpha_2 = 0$ and $\alpha_1\neq 0$,  by choosing $x=\frac{\alpha_1}{\alpha_3}$ we have the representative $\langle \nabla_1 + \nabla_3 \rangle;$

\item if $\lambda = 0$ or $\lambda = 1$, and $\alpha_2 \neq 0$,  by choosing $x=\frac{\alpha_2}{\alpha_3}$ we have the representative $\langle \mu \nabla_1 + \nabla_2 + \nabla_3 \rangle;$

\item if $\lambda\neq 0,1$ and $\alpha_1-\lambda\alpha_2=0$,  by choosing $y=-\frac{x\alpha_2}{(1-\lambda)\lambda \alpha_3}$ and $x=\frac{1}{\sqrt[5]{\alpha_3}}$ we have the representative $\langle \nabla_3 \rangle;$ 

\item if $\lambda\neq 0,1$ and $\alpha_1-\lambda\alpha_2\neq 0$,  by choosing $y=-\frac{x\alpha_2}{(1-\lambda)\lambda \alpha_3}$ and $x=\frac{\alpha_1-\lambda\alpha_2}{\alpha_3}$ we have the representative $\langle \nabla_1 + \nabla_3 \rangle$. 

\end{enumerate}

Therefore, we obtain the families ${\mathcal B}^6_{12}(\lambda)$ and ${\mathcal B}^6_{13}(\lambda)$ associated with the representatives $\langle \nabla_3 \rangle$ and $\langle \nabla_1 + \nabla_3\rangle$, respectively, and the families ${\mathcal B}^6_{14}(\mu)$ and ${\mathcal B}^6_{15}(\mu)$, associated with $\langle \mu \nabla_1 + \nabla_2 + \nabla_3 \rangle$ for the values $\lambda=0$ and $\lambda=1$, respectively.

\begin{longtable}{ll llllll} ${\mathcal B}^6_{12}(\lambda)$ &:&
$e_1 e_1 = e_2$ &   $e_1 e_2=e_3$ & $e_1e_3=e_5$ & $e_1e_5=e_6$ & $e_2 e_1=\lambda e_3 + e_4$  \\
& & $e_2e_2=\lambda e_5$ & $e_2e_3=\lambda e_6$ & $e_3e_1=\lambda e_5$   & $e_3e_2=\lambda e_6$ & $e_5e_1=\lambda e_6$  \\
 ${\mathcal B}^6_{13}(\lambda)$ &:&
$e_1 e_1 = e_2$ &   $e_1 e_2=e_3$ & $e_1e_3=e_5$ & $e_1e_5=e_6$ \\
&& $e_2 e_1=\lambda e_3 + e_4$   & $e_2e_2=\lambda e_5$ & $e_2e_3=\lambda e_6$ & $e_3e_1=\lambda e_5$ \\
&& $e_3e_2=\lambda e_6$ & $e_4e_1=e_6$ 
  & $e_5e_1=\lambda e_6$ \\
 ${\mathcal B}^6_{14}(\mu)$ &:&
$e_1 e_1 = e_2$ &   $e_1 e_2=e_3$ & $e_1e_3=e_5$ & $e_1e_4=e_6$ & $e_1e_5=e_6$  \\
& & $e_2 e_1=e_4$ & $e_2e_2=e_6$ & $e_3e_1=e_6$ & $e_4e_1=\mu e_6$ \\
 ${\mathcal B}^6_{15}(\mu)$ &:&
$e_1 e_1 = e_2$ &   $e_1 e_2=e_3$ & $e_1e_3=e_5$ & $e_1e_4=e_6$ & $e_1e_5=e_6$  \\
&  & $e_2 e_1=e_3+e_4$ & $e_2e_2=e_5+e_6$  & $e_2e_3=e_6$ & $e_3e_1=e_5+e_6$ \\ 
&& $e_3e_2=e_6$  & $e_4e_1=\mu e_6$ & $e_5e_1=e_6$. \\
\end{longtable}

\subsection{Central extensions of ${\mathcal B}^5_{03}(\lambda,1/\lambda)_{\lambda\neq 0}$}

Let us use the following notations:

\[\nabla_1 = [\Delta_{14}], \quad \nabla_2 = [\Delta_{31}], \quad \nabla_3 = \lambda [\Delta_{15}]+[\Delta_{23}]+\lambda[\Delta_{24}]+[\Delta_{32}]+\lambda[\Delta_{42}]+[\Delta_{51}].\]

The automorphism group of ${\mathcal B}^5_{03}(\lambda,1/\lambda)_{\lambda\neq 0}$ consists of invertible matrices of the form
\[\phi=\begin{pmatrix}
x & 0 & 0  & 0 & 0\\
y & x^2 & 0 & 0 & 0\\
z & xy & x^3 & 0 & 0 \\
t & xy & 0 & x^3 & 0\\
s & (1+1/\lambda)xz + (1+\lambda)xt + y^2 & (2+1/\lambda)x^2y & (2+\lambda)x^2y & x^4\\
\end{pmatrix}. \]

Since
\[ \phi^T\begin{pmatrix}
0 & 0 &  0 & \alpha_1 & \lambda\alpha_3\\
0 & 0 & \alpha_3 & \lambda\alpha_3 & 0 \\
\alpha_2 & \alpha_3 & 0 & 0 & 0\\
0 & \lambda\alpha_3 & 0 & 0 & 0\\
\alpha_3 & 0 & 0 & 0 & 0
\end{pmatrix}\phi =
\begin{pmatrix}
\alpha^{*} & \alpha^{**} & \alpha^{***} & \alpha_1^* + \lambda \alpha^{***} & \lambda \alpha_3^* \\
\alpha^{****} & \alpha^{***} & \alpha_3^* & \lambda\alpha_3^* & 0 \\
\alpha_2^*+(1/\lambda) \alpha^{***} & \alpha_3^* & 0 & 0 &0\\
\alpha^{***} & \lambda\alpha_3^* & 0 & 0 &0\\
\alpha_3^* & 0 & 0 & 0 &0\\
\end{pmatrix},
\]
 we have that the action of ${\rm Aut} ({\mathcal B}^5_{03}(\lambda,1/\lambda)_{\lambda\neq 0})$ on the subspace
$\langle \sum\limits_{i=1}^3 \alpha_i\nabla_i  \rangle$
is given by
$\langle \sum\limits_{i=1}^3 \alpha_i^* \nabla_i \rangle,$
where
\[\alpha_1^*=x^4\alpha_1 + (1-\lambda)\lambda x^3y\alpha_3, \quad \alpha_2^* =x^4\alpha_2 + (1-1/\lambda) x^3y\alpha_3, \quad \alpha_3^* = x^5\alpha_3. \]

We can suppose that $\alpha_3\neq0$. We have the following cases:
\begin{enumerate}
\item if $\lambda = 1$, and  $\alpha_1 =\alpha_2= 0$, then we get the representative $\langle \nabla_3 \rangle;$

\item if $\lambda = 1$,  $\alpha_2 = 0$ and $\alpha_1\neq 0$,  by choosing $x=\frac{\alpha_1}{\alpha_3}$ we have the representative $\langle \nabla_1 + \nabla_3 \rangle;$

\item if $\lambda = 1$, and $\alpha_2 \neq 0$,  by choosing $x=\frac{\alpha_2}{\alpha_3}$ we have the representative $\langle \mu \nabla_1 + \nabla_2 + \nabla_3 \rangle;$

\item if $\lambda\neq 1$ and $\alpha_1+\lambda^2\alpha_2=0$,  by choosing $y=-\frac{x\alpha_2}{(1-1/\lambda)\lambda \alpha_3}$ and $x=\frac{1}{\sqrt[5]{\alpha_3}}$ we get the representative $\langle \nabla_3 \rangle;$

\item if $\lambda\neq 1$ and $\alpha_1+\lambda^2\alpha_2\neq 0$,  by choosing $y=-\frac{x\alpha_2}{(1-1/\lambda)\lambda \alpha_3}$ and $x=\frac{\alpha_1+\lambda^2\alpha_2}{\alpha_3}$ we have the representative $\langle \nabla_1 + \nabla_3 \rangle$.

\end{enumerate}

Hence, we get the following new algebras, associated with $\langle \nabla_3 \rangle$ and $\langle \nabla_1 + \nabla_3\rangle$ for every value of $\lambda\neq 0$, and with $\langle \mu \nabla_1 + \nabla_2 + \nabla_3 \rangle$ for $\lambda=1$, respectively:

\begin{longtable}{ll lllllll} ${\mathcal B}^6_{16}(\lambda)_{\lambda\neq 0}$ &:&
$e_1 e_1 = e_2$ &   $e_1 e_2=e_4$ & $e_1e_3=e_5$ & $e_1e_4=\lambda e_5$ & $e_1e_5=\lambda e_6$\\
& & $e_2 e_1=e_3$   & $e_2e_2=e_5$  & $e_2e_3=e_6$ & $e_2e_4=\lambda e_6$ & $e_3e_1=(1/\lambda) e_5$\\
&& $e_3e_2=e_6$ & $e_4e_1=e_5$    & $e_4e_2=\lambda e_6$ & $e_5e_1=e_6$  \\
 ${\mathcal B}^6_{17}(\lambda)_{\lambda\neq 0}$ &:&
$e_1 e_1 = e_2$ &   $e_1 e_2=e_4$ & $e_1e_3=e_5$ & $e_1e_4=\lambda e_5 + e_6$  & $e_1e_5=\lambda e_6$ \\
&&  $e_2 e_1=e_3$   & $e_2e_2=e_5$ & $e_2e_3=e_6$ & $e_2e_4=\lambda e_6$ & $e_3e_1=(1/\lambda) e_5$ \\
&& $e_3e_2=e_6$ & $e_4e_1=e_5$   & $e_4e_2=\lambda e_6$ & $e_5e_1=e_6$ & \\

 ${\mathcal B}^6_{18}(\mu)$ &:&
$e_1 e_1 = e_2$ &   $e_1 e_2=e_4$ & $e_1e_3=e_5$  & $e_1e_4=e_5 + \mu e_6$ & $e_1e_5=e_6$ \\
&& $e_2 e_1=e_3$    & $e_2e_2=e_5$ & $e_2e_3=e_6$ & $e_2e_4=e_6$ & $e_3e_1=e_5+e_6$ \\
&& $e_3e_2=e_6$ & $e_4e_1=e_5$    & $e_4e_2=e_6$ & $e_5e_1=e_6$. & \\
\end{longtable}

\subsection{Central extensions of ${\mathcal B}^5_{05}$}

Let us use the following notations:
\[\nabla_1=[\Delta_{31}], \quad \nabla_2=[\Delta_{14}]+[\Delta_{23}]+[\Delta_{32}]+[\Delta_{51}].\]

The automorphism group of ${\mathcal B}^5_{05}$ consists of invertible matrices of the form
\[\phi=\begin{pmatrix}
1 & 0 & 0  & 0 & 0\\
0 & 1 & 0 & 0 & 0\\
x & 0 & 1 & 0 & 0\\
y & x & 0 & 1 & 0\\
z & x+y & x & 0 & 1\\
\end{pmatrix}. \]

Since
\[ \phi^T\begin{pmatrix}
0 & 0 & 0 & \alpha_2 & 0\\
0 & 0 & \alpha_2 & 0 & 0 \\
\alpha_1 & \alpha_2 & 0 & 0 & 0\\
0 & 0 & 0 & 0 & 0\\
\alpha_2 & 0 & 0 & 0 & 0
\end{pmatrix}\phi =
\begin{pmatrix}
\alpha^{*} & \alpha^{**} & 0 &  \alpha_2^* & 0 \\
\alpha^{***} & 0 &  \alpha_2^* & 0 & 0 \\
\alpha_1^*+\alpha^{**} & \alpha_2^* & 0 & 0 & 0\\
0 & 0 & 0 & 0 & 0\\
\alpha_2^* & 0 & 0 & 0 & 0\\
\end{pmatrix},
\]
 we have that the action of ${\rm Aut} ({\mathcal B}^5_{05})$ on the subspace
$\langle \sum\limits_{i=1}^2 \alpha_i\nabla_i  \rangle$
is given by
$\langle \sum\limits_{i=1}^2 \alpha_i^* \nabla_i \rangle,$
where
\[\alpha_1^* = \alpha_1-x\alpha_2, \quad \alpha_2^* = \alpha_2. \]

We can suppose that $ \alpha_2\neq0$. Then, by choosing $x=\frac{\alpha_1}{\alpha_2}$, we have the representative $\langle \nabla_2 \rangle$, whose associated algebra is the following:

\begin{longtable}{lll lllllll} ${\mathcal B}^6_{19}$ &:&
$e_1 e_1 = e_2$ & $e_1 e_2=e_4$ & $e_1e_3=e_5$ & $e_1 e_4=e_6$ & $e_2 e_1=e_3$ & $e_2 e_2 = e_5$  \\
& & $e_2e_3=e_6$  & $e_3e_1=e_4$  & $e_3e_2=e_6$ & $e_4 e_1=e_5$ & $e_5 e_1 = e_6$. \\
\end{longtable}

\subsection{Central extensions of ${\mathcal B}^5_{06}$}

Let us use the following notations:
\[\nabla_1=[\Delta_{41}], \quad \nabla_2=[\Delta_{13}]+[\Delta_{22}]+[\Delta_{51}].\]

The automorphism group of ${\mathcal B}^5_{06}$ consists of invertible matrices of the form
\[\phi=\begin{pmatrix}
1 & 0 & 0  & 0 & 0\\
x & 1 & 0 & 0 & 0\\
y & x & 1 & 0 & 0\\
z & y & x & 1 & 0\\
t & x+z & y & x & 1\\
\end{pmatrix}. \]

Since
\[ \phi^T\begin{pmatrix}
0 & 0 & \alpha_2 & 0 & 0\\
0 & \alpha_2 & 0 & 0 & 0 \\
0 & 0 & 0 & 0 & 0\\
\alpha_1 & 0 & 0 & 0 & 0\\
\alpha_2 & 0 & 0 & 0 & 0
\end{pmatrix}\phi =
\begin{pmatrix}
\alpha^{*} & \alpha^{**} & \alpha_2^* & 0 & 0 \\
\alpha^{***} & \alpha_2^* & 0 & 0 & 0 \\
\alpha^{****} & 0 & 0 & 0 & 0\\
\alpha_1^*+\alpha^{**} & 0 & 0 & 0 & 0\\
\alpha_2^* & 0 & 0 & 0 & 0\\
\end{pmatrix},
\]
 we have that the action of ${\rm Aut} ({\mathcal B}^5_{06})$ on the subspace
$\langle \sum\limits_{i=1}^2 \alpha_i\nabla_i  \rangle$
is given by
$\langle \sum\limits_{i=1}^2 \alpha_i^* \nabla_i \rangle,$
where
\[\alpha_1^* = \alpha_1-x\alpha_2, \quad \alpha_2^* = \alpha_2. \]

We can suppose that $ \alpha_2\neq0$. Then, by choosing $x=\frac{\alpha_1}{\alpha_2}$, we have the representative $\langle \nabla_2 \rangle$. Hence, we obtain the following new algebra:

\begin{longtable}{ll llllllllll} ${\mathcal B}^6_{20}$ &:&
$e_1 e_1 = e_2$ & $e_1e_2=e_5$ & $e_1 e_3 = e_6$ & $e_2 e_1=e_3$ \\
&& $e_2 e_2=e_6$ & $e_3 e_1=e_4$  & $e_4e_1=e_5$    & $e_5 e_1=e_6$. && \\
\end{longtable}

\subsection{Central extensions of ${\mathcal B}^5_{07}$}

Let us use the following notations:
\[\nabla_1=[\Delta_{12}], \quad \nabla_2=[\Delta_{51}].\]

The automorphism group of ${\mathcal B}^5_{07}$ consists of invertible matrices of the form
\[\phi=\begin{pmatrix}
x & 0 & 0  & 0 & 0\\
y & x^2 & 0 & 0 & 0\\
z & xy & x^3 & 0 & 0\\
t & xz & x^2y & x^4 & 0\\
s & xt & x^2z & x^3y & x^5\\
\end{pmatrix}. \]

Since
\[ \phi^T\begin{pmatrix}
0 & \alpha_1 & 0 & 0 & 0\\
0 & 0 & 0 & 0 & 0 \\
0 & 0 & 0 & 0 & 0\\
0 & 0 & 0 & 0 & 0\\
\alpha_2 & 0 & 0 & 0 & 0
\end{pmatrix}\phi =
\begin{pmatrix}
\alpha^{*} & \alpha_1^* & 0 & 0 & 0 \\
\alpha^{**} & 0 & 0 & 0 & 0 \\
\alpha^{***} & 0 & 0 & 0 & 0\\
\alpha^{****} & 0 & 0 & 0 & 0\\
\alpha_2^* & 0 & 0 & 0 & 0\\
\end{pmatrix},
\]
 we have that the action of ${\rm Aut} ({\mathcal B}^5_{07})$ on the subspace
$\langle \sum\limits_{i=1}^2 \alpha_i\nabla_i  \rangle$
is given by
$\langle \sum\limits_{i=1}^2 \alpha_i^* \nabla_i \rangle,$
where
\[\alpha_1^* = x^3\alpha_1, \quad \alpha_2^* = x^6\alpha_2. \]

We can suppose that $ \alpha_2\neq0$. We have the following cases:
\begin{enumerate}
\item if $\alpha_1 \neq 0$,  by choosing $x=\sqrt[3]{\frac{\alpha_1}{\alpha_2}}$ we have the representative $\langle \nabla_1 + \nabla_2 \rangle;$

\item if $\alpha_1=0$,  by choosing $x=\frac{1}{\sqrt[6]{\alpha_2}}$ we have the representative $\langle \nabla_2 \rangle.$

\end{enumerate}

Hence, we get the following new algebras:

\begin{longtable}{ll llllll} ${\mathcal B}^6_{21}$ &:&
$e_1 e_1 = e_2$ & $e_1e_2=e_6$ & $e_2 e_1=e_3$ & $e_3 e_1=e_4$ & $e_4e_1=e_5$ & $e_5 e_1 = e_6$ \\
 ${\mathcal B}^6_{22}$ &:&
$e_1 e_1 = e_2$ & $e_2 e_1=e_3$ & $e_3 e_1=e_4$ & $e_4e_1=e_5$ & $e_5e_1=e_6$. \\
\end{longtable}

\subsection{Central extensions of ${\mathcal B}^5_{08}$}

Let us use the following notations:
\[\nabla_1=[\Delta_{21}], \quad \nabla_2=[\Delta_{15}]+[\Delta_{23}]+[\Delta_{32}]+[\Delta_{41}].\]

The automorphism group of ${\mathcal B}^5_{08}$ consists of invertible matrices of the form
\[\phi=\begin{pmatrix}
1 & 0 & 0  & 0 & 0\\
0 & 1 & 0 & 0 & 0\\
x & 0 & 1 & 0 & 0\\
y & x & 0 & 1 & 0\\
z & x+y & x & 0 & 1\\
\end{pmatrix}. \]

Since
\[ \phi^T\begin{pmatrix}
0 & 0 & 0 & 0 & \alpha_2\\
\alpha_1 & 0 & \alpha_2 & 0 & 0 \\
0 & \alpha_2 & 0 & 0 & 0\\
\alpha_2 & 0 & 0 & 0 & 0\\
0 & 0 & 0 & 0 & 0
\end{pmatrix}\phi =
\begin{pmatrix}
\alpha^{*} & \alpha^{**} & \alpha^{***} & 0 & \alpha_2^* \\
\alpha_1^*+\alpha^{***} & 0 & \alpha_2^* & 0 & 0 \\
0 & \alpha_2^* & 0 & 0 & 0\\
\alpha_2^* & 0 & 0 & 0 & 0\\
0 & 0 & 0 & 0 & 0\\
\end{pmatrix},
\]
 we have that the action of ${\rm Aut} ({\mathcal B}^5_{08})$ on the subspace
$\langle \sum\limits_{i=1}^2 \alpha_i\nabla_i  \rangle$
is given by
$\langle \sum\limits_{i=1}^2 \alpha_i^* \nabla_i \rangle,$
where
\[\alpha_1^* = \alpha_1+x\alpha_2, \quad \alpha_2^* = \alpha_2. \]

We can suppose that $ \alpha_2\neq0$. Then, by choosing $x=-\frac{\alpha_1}{\alpha_2}$, we have the representative $\langle \nabla_2 \rangle$, with the following associated algebra:

\begin{longtable}{ll llllllll} ${\mathcal B}^6_{23}$ &:&
$e_1e_1=e_2$ & $e_1e_2=e_3$ & $e_1e_3=e_4$  & $e_1e_4=e_5$ & $e_1e_5=e_6$ & $e_2e_1=e_4$ \\
&& $e_2e_2=e_5$ & $e_2e_3=e_6$ & $e_3e_1=e_5$ & $e_3e_2=e_6$ & $e_4e_1=e_6$. \\
\end{longtable}

\subsection{Central extensions of ${\mathcal B}^5_{09}$}

Let us use the following notations:
\[\nabla_1 = [\Delta_{21}], \quad \nabla_2 = [\Delta_{15}] + [\Delta_{23}] + [\Delta_{24}] + [\Delta_{32}]  + [\Delta_{33}] + [\Delta_{41}] + [\Delta_{42}] + [\Delta_{51}].\]

The automorphism group of ${\mathcal B}^5_{09}$ consists of invertible matrices of the form
\[\phi=\begin{pmatrix}
1 & 0 & 0  & 0 & 0\\
0 & 1 & 0 & 0 & 0\\
x & 0 & 1 & 0 & 0\\
y & 2x & 0 & 1 & 0\\
z & x+2y & 3x & 0 & 1\\
\end{pmatrix}. \]

Since
\[ \phi^T\begin{pmatrix}
0 & 0 & 0 & 0 & \alpha_2\\
\alpha_1 & 0 & \alpha_2 & \alpha_2 & 0 \\
0 & \alpha_2 & \alpha_2 & 0 & 0\\
\alpha_2 & \alpha_2 & 0 & 0 & 0\\
\alpha_2 & 0 & 0 & 0 & 0
\end{pmatrix}\phi =
\begin{pmatrix}
\alpha^{*} & \alpha^{**} & \alpha^{***} & 0 & \alpha_2^* \\
\alpha_1^*+\alpha^{**}+\alpha^{***} & \alpha^{***} & \alpha_2^* & \alpha_2^* & 0 \\
\alpha^{***} & \alpha_2^* & \alpha_2^* & 0 & 0\\
\alpha_2^* & \alpha_2^* & 0 & 0 & 0\\
\alpha_2^* & 0 & 0 & 0 & 0\\
\end{pmatrix},
\]
 we have that the action of ${\rm Aut} ({\mathcal B}^5_{09})$ on the subspace
$\langle \sum\limits_{i=1}^2 \alpha_i\nabla_i  \rangle$
is given by
$\langle \sum\limits_{i=1}^2 \alpha_i^* \nabla_i \rangle,$
where
\[\alpha_1^* = \alpha_1-2x\alpha_2, \quad \alpha_2^* = \alpha_2. \]

We can suppose that $ \alpha_2\neq0$. Then, by choosing $x=\frac{\alpha_1}{2\alpha_2}$, we have the representative $\langle \nabla_2 \rangle$. Hence, we get the following new algebra:

\begin{longtable}{ll llllllll} ${\mathcal B}^6_{24}$ &:&
$e_1e_1=e_2$ & $e_1e_2=e_3$ & $e_1e_3=e_4$  & $e_1e_4=e_5$ & $e_1e_5=e_6$ \\
&& $e_2e_1=e_3+e_4$  & $e_2e_2=e_4+e_5$ & $e_2e_3=e_5+e_6$ & $e_2e_4=e_6$ & $e_3e_1=e_4+e_5$ \\
&& $e_3e_2=e_5+e_6$  & $e_3e_3=e_6$   &  $e_4e_1=e_5+e_6$ & $e_4e_2=e_6$ & $e_5e_1=e_6$. \\
\end{longtable}

\subsection{Central extensions of ${\mathcal B}^5_{10}(\lambda)$}

Let us use the following notations:
\[\nabla_1 = [\Delta_{21}], \quad \nabla_2 = [\Delta_{15}] + \lambda[\Delta_{24}] + \lambda[\Delta_{33}] + \lambda[\Delta_{42}] + \lambda[\Delta_{51}].\]

The automorphism group of ${\mathcal B}^5_{10}(\lambda)$ consists of invertible matrices of the form
\[\phi=\begin{pmatrix}
x & 0 & 0  & 0 & 0\\
0 & x^2 & 0 & 0 & 0\\
0 & 0 & x^3 & 0 & 0\\
y & 0 & 0 & x^4 & 0\\
z & (1+\lambda)xy & 0 & 0 & x^5\\
\end{pmatrix}. \]

Since
\[ \phi^T\begin{pmatrix}
0 & 0 & 0 & 0 & \alpha_2\\
\alpha_1 & 0 & 0 & \lambda\alpha_2 & 0 \\
0 & 0 & \lambda\alpha_2 & 0 & 0\\
0 & \lambda\alpha_2 & 0 & 0 & 0\\
\lambda\alpha_2 & 0 & 0 & 0 & 0
\end{pmatrix}\phi =
\begin{pmatrix}
\alpha^{*} & \alpha^{**} & 0 & 0 & \alpha_2^* \\
\alpha_1^*+\lambda \alpha^{**} & 0 & 0 & \lambda\alpha_2^* & 0 \\
0 & 0 & \lambda\alpha_2^* & 0 & 0\\
0 & \lambda\alpha_2^* & 0 & 0 & 0\\
\lambda\alpha_2^* & 0 & 0 & 0 & 0\\
\end{pmatrix},
\]
 we have that the action of ${\rm Aut} ({\mathcal B}^5_{10}(\lambda))$ on the subspace
$\langle \sum\limits_{i=1}^2 \alpha_i\nabla_i  \rangle$
is given by
$\langle \sum\limits_{i=1}^2 \alpha_i^* \nabla_i \rangle,$
where
\[\alpha_1^* = x^3\alpha_1+(1-\lambda)\lambda x^2y\alpha_2, \quad \alpha_2^* = x^6\alpha_2. \]

We can suppose that $ \alpha_2\neq0$. We have the following cases:
\begin{enumerate}
\item if $\lambda=0$ or $\lambda = 1$, and $\alpha_1 \neq 0$,  by choosing $x=\sqrt[3]{\frac{\alpha_1}{\alpha_2}}$ we have the representative $\langle \nabla_1 + \nabla_2 \rangle$;

\item if $\lambda=0$ or $\lambda = 1$, and $\alpha_1 = 0$,  by choosing $x=\frac{1}{\sqrt[6]{\alpha_2}}$ we have the representative $\langle \nabla_2 \rangle$;

\item if $\lambda \neq 0, 1$,  by choosing $y=-\frac{x\alpha_1}{(1-\lambda)\lambda\alpha_2}$ and $x=\frac{1}{\sqrt[6]{\alpha_2}}$ we have the representative $\langle \nabla_2 \rangle.$

\end{enumerate}

Hence, we obtain the following new algebras, associated with $\langle \nabla_2 \rangle$ for every value of $\lambda$, and with $\nabla_1 + \nabla_2$ for $\lambda=0$ or $\lambda=1$, respectively.

\begin{longtable}{ll llllllll}  ${\mathcal B}^6_{25}(\lambda)$ &:&
$e_1e_1=e_2$ & $e_1e_2=e_3$ & $e_1e_3=e_4$ & $e_1e_4=e_5$ & $e_1e_5=e_6$ & $e_2e_1=\lambda e_3$ \\
& & $e_2e_2=\lambda e_4$ & $e_2e_3=\lambda e_5$ & $e_2e_4=\lambda e_6$ & $e_3e_1=\lambda e_4$ & $e_3e_2=\lambda e_5$  & $e_3e_3=\lambda e_6$ \\
&  & $e_4e_1=\lambda e_5$ & $e_4e_2=\lambda e_6$ & $e_5e_1=\lambda e_6$.\\
${\mathcal B}^6_{26}$ &:&
$e_1e_1=e_2$ & $e_1e_2=e_3$ & $e_1e_3=e_4$  & $e_1e_4=e_5$ & $e_1e_5=e_6$ & $e_2e_1=e_6$ \\
 ${\mathcal B}^6_{27}$ &:&
$e_1e_1=e_2$ & $e_1e_2=e_3$ & $e_1e_3=e_4$ & $e_1e_4=e_5$ & $e_1e_5=e_6$ & $e_2e_1=e_3+e_6$ \\
& & $e_2e_2=e_4$ & $e_2e_3=e_5$  & $e_2e_4=e_6$ & $e_3e_1=e_4$ & $e_3e_2=e_5$ & $e_3e_3=e_6$ \\
& & $e_4e_1=e_5$ & $e_4e_2=e_6$ & $e_5e_1=e_6$ \\
\end{longtable}

\subsection{Central extensions of ${\mathcal B}^5_{11}$}

Let us use the following notations:
\[\nabla_1 = [\Delta_{21}], \quad \nabla_2 = [\Delta_{15}]+[\Delta_{22}]+[\Delta_{31}].\]

The automorphism group of ${\mathcal B}^5_{11}$ consists of invertible matrices of the form
\[\phi=\begin{pmatrix}
1 & 0 & 0  & 0 & 0\\
x & 1 & 0 & 0 & 0\\
y & x & 1 & 0 & 0\\
z & y & x & 1 & 0\\
t & x+z & y & x & 1\\
\end{pmatrix}. \]

Since
\[ \phi^T\begin{pmatrix}
0 & 0 & 0 & 0 & \alpha_2\\
\alpha_1 & \alpha_2 & 0 & 0 & 0 \\
\alpha_2 & 0 & 0 & 0 & 0\\
0 & 0 & 0 & 0 & 0\\
0 & 0 & 0 & 0 & 0
\end{pmatrix}\phi =
\begin{pmatrix}
\alpha^{*} & \alpha^{**} & \alpha^{***} & \alpha^{****} & \alpha_2^* \\
\alpha_1^*+\alpha^{****} & \alpha_2^* & 0 & 0 & 0 \\
\alpha_2^* & 0 & 0 & 0 & 0\\
0 & 0 & 0 & 0 & 0\\
0 & 0 & 0 & 0 & 0\\
\end{pmatrix},
\]
 we have that the action of ${\rm Aut} ({\mathcal B}^5_{10}(\lambda))$ on the subspace
$\langle \sum\limits_{i=1}^2 \alpha_i\nabla_i  \rangle$
is given by
$\langle \sum\limits_{i=1}^2 \alpha_i^* \nabla_i \rangle,$
where
\[\alpha_1^* = \alpha_1+x\alpha_2, \quad \alpha_2^* = \alpha_2. \]

We can suppose that $ \alpha_2\neq0$. Then, by taking $x=-\frac{\alpha_1}{\alpha_2}$, we have the representative $\langle \nabla_2 \rangle$, with the following associated algebra:

\begin{longtable}{ll llllllll} ${\mathcal B}^6_{28}$ &:&
$e_1e_1=e_2$ & $e_1e_2=e_3$ & $e_1e_3=e_4$  & $e_1e_4=e_5$ \\
&& $e_1e_5=e_6$ & $e_2e_1=e_5$ & $e_2e_2=e_6$  & $e_3e_1=e_6$. &&  \\
\end{longtable}

\subsection{Central extensions of ${\mathcal B}^5_{12}$}

Let us use the following notations:
\[\nabla_1 = [\Delta_{21}], \quad \nabla_2 =[\Delta_{15}]+[\Delta_{22}]+[\Delta_{24}]+[\Delta_{31}]+[\Delta_{33}]+[\Delta_{42}]+[\Delta_{51}].\]

The automorphism group of ${\mathcal B}^5_{12}$ consists of invertible matrices of the form
\[\phi=\begin{pmatrix}
1 & 0 & 0  & 0 & 0\\
x & 1 & 0 & 0 & 0\\
y & 2x & 1 & 0 & 0\\
z & x^2+2y & 3x & 1 & 0\\
t & x(1+2y)+2z & 3x^2+3y & 4x & 1\\
\end{pmatrix}. \]

Since
\[ \phi^T\begin{pmatrix}
0 & 0 & 0 & 0 & \alpha_2\\
\alpha_1 & \alpha_2 & 0 & \alpha_2 & 0 \\
\alpha_2 & 0 & \alpha_2 & 0 & 0\\
0 & \alpha_2 & 0 & 0 & 0\\
\alpha_2 & 0 & 0 & 0 & 0
\end{pmatrix}\phi =
\begin{pmatrix}
\alpha^{*} & \alpha^{**} & \alpha^{***} & \alpha^{****} & \alpha_2^* \\
\alpha_1^*+\alpha^{**}+\alpha^{****} & \alpha_2^*+\alpha^{***} & \alpha^{****} & \alpha_2^* & 0 \\
\alpha_2^*+\alpha^{***} & \alpha^{****} & \alpha_2^* & 0 & 0\\
\alpha^{****} & \alpha_2^* & 0 & 0 & 0\\
\alpha_2^* & 0 & 0 & 0 & 0\\
\end{pmatrix},
\]
 we have that the action of ${\rm Aut} ({\mathcal B}^5_{10}(\lambda))$ on the subspace
$\langle \sum\limits_{i=1}^2 \alpha_i\nabla_i  \rangle$
is given by
$\langle \sum\limits_{i=1}^2 \alpha_i^* \nabla_i \rangle,$
where
\[\alpha_1^* = \alpha_1-3x\alpha_2, \quad \alpha_2^* = \alpha_2. \]

We can suppose that $ \alpha_2\neq0$. Then, by taking $x=\frac{\alpha_1}{3\alpha_2}$, we have the representative $\langle \nabla_2 \rangle$.
Hence, we obtain the following new algebra:

\begin{longtable}{ll lllllll} ${\mathcal B}^6_{29}$ &:&
$e_1e_1=e_2$ & $e_1e_2=e_3$ & $e_1e_3=e_4$  & $e_1e_4=e_5$ & $e_1e_5=e_6$ \\
&& $e_2e_1=e_3+e_5$   & $e_2e_2=e_4+e_6$ & $e_2e_3=e_5$ & $e_2e_4=e_6$ & $e_3e_1=e_4+e_6$ \\
&& $e_3e_2=e_5$ & $e_3e_3=e_6$  & $e_4e_1=e_5$ & $e_4e_2=e_6$ & $e_5e_1=e_6$. \\
\end{longtable}

\subsection{Classification theorem}
Summarizing the results of the present section, we obtain the following theorem.

\begin{theoremB}
Let $\mathcal B$ be a $6$-dimensional complex one-generated nilpotent bicommutative algebra.
Then, $\mathcal B$ is isomorphic to an algebra from the following list:
 
\begin{longtable}{l|lllll}
 
\hline ${\mathcal B}^6_{01}$ &
$e_1 e_1 = e_2$ &    $e_1e_2=e_4$ & $e_1 e_4=e_5$ & $e_2 e_1=e_3$ & $e_3 e_1=e_6$ \\

\hline ${\mathcal B}^6_{02}(\lambda)$ &
$e_1 e_1 = e_2$ &   $e_1 e_2=e_4$ & $e_1e_3=e_6$ &  $e_1 e_4=e_5$ \\
&  $e_2 e_1= e_3$  
& $e_2e_2= e_6$ & $e_3e_1=\lambda e_6$ &  $e_4e_1= e_6$    \\

\hline ${\mathcal B}^6_{03}(\lambda, \mu)$ & 
$e_1 e_1 = e_2$ &   $e_1 e_2=e_4$ &  $e_1e_3=  e_6$ &  \multicolumn{2}{l }{ $e_1 e_4=\lambda e_5 + \mu e_6$} \\ 
& $e_2 e_1= e_3$ &  $e_2e_2= e_6$ & $e_3e_1= e_5$ & $e_4e_1= e_6$  &  \\

\hline ${\mathcal B}^6_{04}$ &
$e_1 e_1 = e_2$ &    $e_1 e_2=e_4$ & $e_1 e_3=e_6$ & $e_2e_1=e_3$ \\& $e_2e_2=e_6$ 
&   $e_3e_1=e_4+e_5$ &  $e_4 e_1=e_6$ \\

\hline ${\mathcal B}^6_{05}$ &
$e_1 e_1 = e_2$ & $e_1e_2=e_5$ &   $e_2 e_1=e_3$ &    $e_3 e_1=e_4$  & $e_4 e_1=e_6$  \\ 

\hline ${\mathcal B}^6_{06}$ & 
$e_1 e_1 = e_2$ &  $e_1 e_2=e_3$ & $e_1 e_3=e_4$ & $e_1e_4=e_6$ \\& $e_2e_1=e_4+e_5$  
&  $e_2e_2=e_6$ & $e_3e_1=e_6$ \\

\hline ${\mathcal B}^6_{07}$ &
$e_1 e_1 = e_2$ &   $e_1 e_2=e_3$ & $e_1 e_3=e_4$  & $e_1e_4=e_6$  
& $e_2e_1= e_3+e_4+e_5$  \\
& $e_2e_2=e_4+e_6$ & $e_2e_3=e_6$  & $e_3e_1=e_4+e_6$  & $e_3e_2=e_6$ & $e_4e_1=e_6$\\

\hline ${\mathcal B}^6_{08}(\lambda)$ & 
$e_1 e_1 = e_2$ &   $e_1 e_2=e_3$ & $e_1 e_3=e_4$ & $e_1e_4=e_6$ & $e_2e_1=\lambda e_3+e_5$  \\ 
& $e_2e_2=\lambda e_4$ & $e_2e_3=\lambda e_6$ &  $e_3e_1=\lambda e_4$  & $e_3e_2=\lambda e_6$ & $e_4e_1=\lambda e_6$ \\

\hline ${\mathcal B}^6_{09}(\lambda)$ & 
$e_1 e_1 = e_2$  & $e_1 e_2=e_4$ & $e_1e_3=e_6$ & $e_1e_4=\lambda e_6$ & $e_2 e_1=e_3$  \\ 
& $e_2e_2=e_6$ & $e_3e_1=e_5$ & $e_4e_1=e_6$ & $e_5e_1=e_6$ & \\

\hline ${\mathcal B}^6_{10}$ & 
$e_1 e_1 = e_2$ & $e_1 e_2=e_4$ & $e_1e_4=e_6$ \\ & $e_2 e_1=e_3$ & $e_3e_1=e_5$    
& $e_5e_1=e_6$ &   \\

\hline ${\mathcal B}^6_{11}$ & 
$e_1 e_1 = e_2$  & $e_1 e_2=e_4$ & $e_2 e_1=e_3$ & $e_3e_1=e_5$ & $e_5e_1=e_6$ \\ 

\hline ${\mathcal B}^6_{12}(\lambda)$ &
$e_1 e_1 = e_2$ &   $e_1 e_2=e_3$ & $e_1e_3=e_5$ & $e_1e_5=e_6$ & $e_2 e_1=\lambda e_3 + e_4$ \\
& $e_2e_2=\lambda e_5$ & $e_2e_3=\lambda e_6$ & $e_3e_1=\lambda e_5$   & $e_3e_2=\lambda e_6$ & $e_5e_1=\lambda e_6$ \\

\hline ${\mathcal B}^6_{13}(\lambda)$ &
$e_1 e_1 = e_2$ &   $e_1 e_2=e_3$ & $e_1e_3=e_5$ & $e_1e_5=e_6$ \\
& $e_2 e_1=\lambda e_3 + e_4$ 
& $e_2e_2=\lambda e_5$ & $e_2e_3=\lambda e_6$ & $e_3e_1=\lambda e_5$ \\
& $e_3e_2=\lambda e_6$ & $e_4e_1=e_6$ & $e_5e_1=\lambda e_6$   \\

\hline ${\mathcal B}^6_{14}(\lambda)$ &
$e_1 e_1 = e_2$ &   $e_1 e_2=e_3$ & $e_1e_3=e_5$ & $e_1e_4=e_6$ & $e_1e_5=e_6$  \\
& $e_2 e_1=e_4$ & $e_2e_2=e_6$ & $e_3e_1=e_6$ & $e_4e_1=\lambda e_6$ & \\

\hline ${\mathcal B}^6_{15}(\lambda)$ &
$e_1 e_1 = e_2$ &   $e_1 e_2=e_3$ & $e_1e_3=e_5$ & $e_1e_4=e_6$ \\
& $e_1e_5=e_6$  & $e_2 e_1=e_3+e_4$ & $e_2e_2=e_5+e_6$  & $e_2e_3=e_6$ \\
& $e_3e_1=e_5+e_6$ & $e_3e_2=e_6$ & $e_4e_1=\lambda e_6$ & $e_5e_1=e_6$  \\

\hline ${\mathcal B}^6_{16}(\lambda)_{\lambda\neq 0}$ &
$e_1 e_1 = e_2$ &   $e_1 e_2=e_4$ & $e_1e_3=e_5$ & $e_1e_4=\lambda e_5$ & $e_1e_5=\lambda e_6$ \\ 
& $e_2 e_1=e_3$ & $e_2e_2=e_5$  & $e_2e_3=e_6$ & $e_2e_4=\lambda e_6$ & $e_3e_1=(1/\lambda) e_5$\\
 & $e_3e_2=e_6$ & $e_4e_1=e_5$ & $e_4e_2=\lambda e_6$ & $e_5e_1=e_6$ & \\

\hline ${\mathcal B}^6_{17}(\lambda)_{\lambda\neq 0}$ &
$e_1 e_1 = e_2$ &   $e_1 e_2=e_4$ & $e_1e_3=e_5$ & \multicolumn{2}{l}{$e_1e_4=\lambda e_5 + e_6$}\\ 
& $e_1e_5=\lambda e_6$ & $e_2 e_1=e_3$ & $e_2e_2=e_5$ & $e_2e_3=e_6$ & $e_2e_4=\lambda e_6$ \\& $e_3e_1=(1/\lambda) e_5$ 
& $e_3e_2=e_6$ & $e_4e_1=e_5$  & $e_4e_2=\lambda e_6$ & $e_5e_1=e_6$   \\

\hline ${\mathcal B}^6_{18}(\lambda)$ &
$e_1 e_1 = e_2$ &   $e_1 e_2=e_4$ & $e_1e_3=e_5$  & \multicolumn{2}{l}{$e_1e_4=e_5 + \lambda e_6$}\\ 
& $e_1e_5=e_6$  
& $e_2 e_1=e_3$  & $e_2e_2=e_5$ & $e_2e_3=e_6$ & $e_2e_4=e_6$ \\
& $e_3e_1=e_5+e_6$ 
 & $e_3e_2=e_6$ & $e_4e_1=e_5$  & $e_4e_2=e_6$ & $e_5e_1=e_6$  \\

\hline ${\mathcal B}^6_{19}$ &
$e_1 e_1 = e_2$ & $e_1 e_2=e_4$ & $e_1e_3=e_5$ & $e_1 e_4=e_6$ \\
& $e_2 e_1=e_3$  
& $e_2 e_2 = e_5$   & $e_2e_3=e_6$ & $e_3e_1=e_4$ \\
& $e_3e_2=e_6$ & $e_4 e_1=e_5$ & $e_5 e_1 = e_6$  \\

\hline ${\mathcal B}^6_{20}$ &
$e_1 e_1 = e_2$ & $e_1e_2=e_5$ & $e_1 e_3 = e_6$ & $e_2 e_1=e_3$ \\
& $e_2 e_2=e_6$ & $e_3 e_1=e_4$  & $e_4e_1=e_5$  & $e_5 e_1=e_6$  \\

\hline ${\mathcal B}^6_{21}$ &
$e_1 e_1 = e_2$ & $e_1e_2=e_6$ & $e_2 e_1=e_3$ \\
& $e_3 e_1=e_4$ & $e_4e_1=e_5$  & $e_5 e_1 = e_6$  \\

\hline ${\mathcal B}^6_{22}$ &
$e_1 e_1 = e_2$ & $e_2 e_1=e_3$ & $e_3 e_1=e_4$ & $e_4e_1=e_5$ & $e_5e_1=e_6$ \\

\hline ${\mathcal B}^6_{23}$ &
$e_1e_1=e_2$ & $e_1e_2=e_3$ & $e_1e_3=e_4$  & $e_1e_4=e_5$ \\
& $e_1e_5=e_6$ & $e_2e_1=e_4$ & $e_2e_2=e_5$ & $e_2e_3=e_6$ \\
& $e_3e_1=e_5$ & $e_3e_2=e_6$ & $e_4e_1=e_6$ \\

\hline ${\mathcal B}^6_{24}$ &
$e_1e_1=e_2$ & $e_1e_2=e_3$ & $e_1e_3=e_4$  & $e_1e_4=e_5$ & $e_1e_5=e_6$  \\
& $e_2e_1=e_3+e_4$ & $e_2e_2=e_4+e_5$ & $e_2e_3=e_5+e_6$ & $e_2e_4=e_6$ & $e_3e_1=e_4+e_5$  \\
& $e_3e_2=e_5+e_6$  & $e_3e_3=e_6$ &  $e_4e_1=e_5+e_6$ & $e_4e_2=e_6$ & $e_5e_1=e_6$ \\

\hline ${\mathcal B}^6_{25}(\lambda)$ &
$e_1e_1=e_2$ & $e_1e_2=e_3$ & $e_1e_3=e_4$ & $e_1e_4=e_5$ & $e_1e_5=e_6$ \\
& $e_2e_1=\lambda e_3$ & $e_2e_2=\lambda e_4$ & $e_2e_3=\lambda e_5$ & $e_2e_4=\lambda e_6$ & $e_3e_1=\lambda e_4$  \\
& $e_3e_2=\lambda e_5$  & $e_3e_3=\lambda e_6$ & $e_4e_1=\lambda e_5$ & $e_4e_2=\lambda e_6$ & $e_5e_1=\lambda e_6$ \\

\hline ${\mathcal B}^6_{26}$ &
$e_1e_1=e_2$ & $e_1e_2=e_3$ & $e_1e_3=e_4$  \\
& $e_1e_4=e_5$ & $e_1e_5=e_6$ & $e_2e_1=e_6$ & \\

\hline ${\mathcal B}^6_{27}$ &
$e_1e_1=e_2$ & $e_1e_2=e_3$ & $e_1e_3=e_4$ & $e_1e_4=e_5$ & $e_1e_5=e_6$  \\
& $e_2e_1=e_3+e_6$ & $e_2e_2=e_4$ & $e_2e_3=e_5$  & $e_2e_4=e_6$ & $e_3e_1=e_4$ \\
& $e_3e_2=e_5$ & $e_3e_3=e_6$ & $e_4e_1=e_5$ & $e_4e_2=e_6$ & $e_5e_1=e_6$ \\

\hline ${\mathcal B}^6_{28}$ &
$e_1e_1=e_2$ & $e_1e_2=e_3$ & $e_1e_3=e_4$  & $e_1e_4=e_5$ \\
& $e_1e_5=e_6$ & $e_2e_1=e_5$ & $e_2e_2=e_6$ & $e_3e_1=e_6$  \\

\hline ${\mathcal B}^6_{29}$ &
$e_1e_1=e_2$ & $e_1e_2=e_3$ & $e_1e_3=e_4$  & $e_1e_4=e_5$ & $e_1e_5=e_6$ \\
& $e_2e_1=e_3+e_5$ & $e_2e_2=e_4+e_6$ & $e_2e_3=e_5$ & $e_2e_4=e_6$ & $e_3e_1=e_4+e_6$ \\
& $e_3e_2=e_5$ & $e_3e_3=e_6$ & $e_4e_1=e_5$ & $e_4e_2=e_6$ & $e_5e_1=e_6$ \\

\hline

\end{longtable} 
\end{theoremB}


\begin{thebibliography}{99}




\bibitem{ack}
Abdelwahab H.,  Calder\'on A.J., Kaygorodov I.,
    The algebraic and geometric classification of nilpotent binary Lie algebras, 
    International Journal of Algebra and Computation,  29 (2019), 6, 1113--1129.
 
 
\bibitem{kkk18}
      Adashev J.,   Kaygorodov I.,   Khudoyberdiyev A.,   Sattarov A.,  
The algebraic and geometric classification of nilpotent right commutative algebras, 
 Results in Mathematics,   76   (2021), 1, Paper: 24.
 


\bibitem{DKS09} Burde D., Dekimpe K., Deschamps S.,
LR-algebras,  
    New Developments in Lie Theory and Geometry, Amer. Math. Soc., Providence, RI, Contemporary Mathematics 491 (2009), 125--140.


\bibitem{DKV10} Burde D., Dekimpe K., Vercammen K.,
    Complete LR-structures on solvable Lie algebras, 
    Journal of Group Theory, 13 (2010),  5, 703--719.

  


\bibitem{cayley} Cayley A.,
On the theory of analytical forms called trees,
    Phil. Mag., 13 (1857), 19--30;
    Collected Math. Papers, University Press, Cambridge, 3 (1890), 242--246.

\bibitem{cfk182}
Calderón Martín A., Fern\'andez Ouaridi A., Kaygorodov I.,
    The classification of $n$-dimensional anticommutative algebras with $(n-3)$-dimensional annihilator, 
    Communications in Algebra, 47 (2019), 1, 173--181.



\bibitem{cfk19}
Calderón Martín A.,  Fern\'andez Ouaridi A., Kaygorodov I.,
    The classification of $2$-dimensional rigid algebras,
    Linear and Multilinear Algebra,   68  (2020),  4, 828--844.

\bibitem{cfk18}
Calderón Martín A., Fern\'andez Ouaridi A., Kaygorodov I.,
  On the classification of bilinear maps with   radical of a fixed codimension,  
    Linear and Multilinear Algebra, 2020, DOI:10.1080/03081087.2020.1849001

   
\bibitem{kkl18}
  Camacho L., Karimjanov I., Kaygorodov I., Khudoyberdiyev  A.,  
    Central extensions of filiform Zinbiel algebras, 
    Linear and Multilinear Algebra, 2020, DOI: 10.1080/03081087.2020.1764903
 
\bibitem{ckkk19}
Camacho L.,  Karimjanov I.,  Kaygorodov I., Khudoyberdiyev A.,
    One-generated nilpotent Novikov algebras, 
    Linear and Multilinear Algebra,   70  (2022), 2,  331--365.

\bibitem{lisa}
Camacho L., Kaygorodov I.,  Lopatkin V., Salim M., 
    The variety of dual Mock-Lie algebras, Communications in Mathematics,    28  (2020), 2,  161--178. 
 

\bibitem{degr3}
Cicalò S., De Graaf W.,   Schneider C.,
 Six-dimensional nilpotent Lie algebras,
 Linear Algebra and its Applications, 436 (2012), 1, 163--189.


\bibitem{usefi1}
Darijani I., Usefi H.,
 The classification of 5-dimensional $p$-nilpotent restricted Lie algebras over perfect fields, I.,
 Journal of Algebra, 464 (2016), 97--140.


\bibitem{degr2}
De Graaf W., 
 Classification of 6-dimensional nilpotent Lie algebras over fields of characteristic not $2$, 
 Journal of Algebra, 309  (2007), 2, 640--653.

\bibitem{degr1}
De Graaf W., 
 Classification of nilpotent associative algebras of small dimension,
 International Journal of Algebra and Computation, 28 (2018),  1, 133--161.

 \bibitem{karel}
Dekimpe K.,  Ongenae V.,
    Filiform left-symmetric algebras,
    Geometriae Dedicata, 74 (1999),  2, 165--199.


\bibitem{drensky1}
Drensky V., Zhakhayev B.,
 Noetherianity and Specht problem for varieties of bicommutative algebras,
 Journal of Algebra, 499 (2018), 1, 570--582.

\bibitem{drensky2}
Drensky V.,
 Varieties of bicommutative algebras,
 Serdica Mathematical Journal, 45 (2019), 167--188
    
\bibitem{dt03}
Dzhumadildaev A., Tulenbaev K.,
 Bicommutative algebras,
 Russian Mathematical Surveys, 58 (2003), 6, 1196--1197.

\bibitem{dit11}
Dzhumadildaev A., Ismailov N., Tulenbaev K.,
 Free bicommutative algebras,
 Serdica Mathematical Journal, 37 (2011),  1, 25--44.

\bibitem{di}
Dzhumadildaev A., Ismailov N.,
 Polynomial identities of bicommutative algebras, Lie and Jordan elements,
 Communications in Algebra, 46 (2018), 12,  5241--5251.

\bibitem{fkkv19}
Fern\'andez Ouaridi A.,  Kaygorodov I.,  Khrypchenko M., Volkov Yu., 
    Degenerations of nilpotent algebras,
    Journal of Pure and Applied Algebra,   226  (2022),  3,  Paper: 106850.
 

  

\bibitem{ha16}
Hegazi A., Abdelwahab H.,
    Classification of five-dimensional nilpotent Jordan algebras,
    Linear Algebra and its Applications, 494 (2016), 165--218.




\bibitem{hac16}
Hegazi A., Abdelwahab H., Calderón Martín A.,
    The classification of $n$-dimensional non-Lie Malcev algebras with $(n-4)$-dimensional annihilator, 
    Linear Algebra and its Applications, 505 (2016), 32--56.

\bibitem{hac18}
Hegazi A., Abdelwahab H.,  Calderón Martín A.,
    Classification of nilpotent Malcev algebras of small dimensions over arbitrary fields of characteristic not $2$,
    Algebras and Representation Theory, 21 (2018), 1, 19--45.

\bibitem{ikm19}
 Ismailov N.,  Kaygorodov I.,  Mashurov F.,
    The algebraic and geometric classification of nilpotent assosymmetric algebras,
    Algebras and Representation Theory,   24  (2021), 1, 135--148. 

\bibitem{jkk19}
Jumaniyozov D.,  Kaygorodov I.,   Khudoyberdiyev A., 
      The algebraic and geometric classification of nilpotent noncommutative Jordan  algebras,
     Journal of Algebra and its Applications,   20  (2021),  11, 2150202.


\bibitem{gkk}
   Jumaniyozov D., Kaygorodov I.,   Khudoyberdiyev  A.,   
The algebraic  classification of nilpotent commutative $\mathfrak{CD}$-algebras,   Communications in Algebra,   49  (2021), 4,  1464--1494.
  

\bibitem{krs19}
 Kaygorodov I., Khrypchenko M.,  
The algebraic  classification of nilpotent $\mathfrak{CD}$-algebras, 
 Linear and Multilinear Algebra, 2020, 10.1080/03081087.2020.1856030 
 

\bibitem{kkl19}
Kaygorodov I., Khrypchenko M., Lopes S.,  
    The algebraic and geometric classification of nilpotent anticommutative  algebras,
    Journal of Pure and Applied Algebra, 224 (2020), 8, 106337.
 
\bibitem{ha17}
Kaygorodov I., Khrypchenko M., Lopes S., 
The algebraic      classification of nilpotent  algebras,
 Journal of Algebra and its Applications, 2021,
 DOI: 10.1142/S0219498823500093

\bibitem{kkp19geo}
Kaygorodov I.,  Khrypchenko M.,  Popov Yu.,
    The algebraic and geometric classification of nilpotent terminal algebras,
   Journal of Pure and Applied Algebra,  225 (2021), 6, 106625.

   
\bibitem{kks19}
Kaygorodov I., Khudoyberdiyev A., Sattarov A.,
    One-generated nilpotent terminal algebras, 
       Communications in Algebra,   48  (2020),  10, 4355--4390. 


\bibitem{omirov}
 Kaygorodov I., Lopes S., P\'{a}ez-Guill\'{a}n P., 
Non-associative central extensions of null-filiform associative algebras,
Journal of Algebra,   560 (2020),   1190--1210.

\bibitem{gkks}
Kaygorodov I., Mashurov F., 
One-generated nilpotent assosymmetric algebras, 
Journal of Algebra and its Applications,  21  (2022),  2, 2250031, 21 pp.

\bibitem{kpv19}
Kaygorodov I.,  P\'{a}ez-Guill\'{a}n P.,  Voronin V.,
    The algebraic and geometric classification of nilpotent bicommutative algebras,
    Algebras and Representation Theory,   23  (2020), 6, 2331-2347.
   
 
 
  
\bibitem{kv16}
Kaygorodov I.,   Volkov Yu.,
    The variety of $2$-dimensional algebras over an algebraically closed field,
    Canadian  Journal of Mathematics,  71 (2019),  4, 819--842.
 
\bibitem{bakhrom}	
Masutova K., Omirov B., 
    On some zero-filiform algebras, 
    Ukrainian Mathematical Journal, 66 (2014), 4, 541--552.

\bibitem{ss78}
Skjelbred T., Sund T.,
    Sur la classification des algebres de Lie nilpotentes,
    C. R. Acad. Sci. Paris Ser. A-B, 286 (1978), 5,  A241--A242.

\bibitem{zusmanovich}
Zusmanovich P., 
    Central extensions of current algebras,
    Transactions of the American Mathematical Society, 334 (1992),  1, 143--152.


\end{thebibliography}
\end{document}